\documentstyle[12pt,righttag,amscd,amssymb,a4]{amsart}
\input{xy}
\xyoption{all}
\catcode `\@=11
\def\numberbysection{\@addtoreset{equation}{section}
         \renewcommand{\theequation}{\thesection.\arabic{equation}}}
\numberbysection

\evensidemargin \oddsidemargin
\def\subsubsection{\@startsection{subsubsection}{3}%
  \normalparindent{.5\linespacing\@plus.7\linespacing}{-.5em}%
  {\normalfont\bfseries}}

\def\bitm{\bibitem}

\def\mbn{\overline{M}_{0n}}
\def\mbnn{\overline{M}_{0n-1}}
\def\mbs{\overline{M}_{0S}}

\def\la{\lambda}

\def\del{\partial}

\def\Mb{\overline{M}}

\def\n {\{1,\ldots ,n \}}

\def\med{\smallskip}
\def\sm{\smallskip}

\def\ra{\rightarrow}
\def\shave{}
\def\eps{\epsilon}

\def\beq{\begin{equation}}
\def\eeq{\end{equation}}

\def\bsp{\begin {split}}
\def\espl{\end{split}}
\def\bal{\begin{align}}
\def\eal{\end{align}}
\def\bald{\begin{aligned}}
\def\eald{\end{aligned}}
\def\nn{\nonumber}
\def\Cal{\cal}

\def\A{{\Cal A}}
\def\NIC{\,\text{\bf CohFT}_{\bold1}(k)\,}

\def\LR{\;\longleftrightarrow\;}

\def\cal{\mathcal}
\begin{document}

\title{The tensor product in the theory of Frobenius manifolds} 

\author[Ralph M. Kaufmann]{Ralph M. Kaufmann\\
Max--Planck--Institut f\"ur Mathematik, Bonn}

\address{Max--Planck--Institut f\"ur Mathematik, Bonn}

\begin{abstract}
We introduce the operation of forming the
tensor product in the theory of analytic Frobenius manifolds. 
Building on the results for formal Frobenius manifolds which we
extend to the additional structures of Euler fields and flat identities, 
we prove that
the tensor product of pointed germs of Frobenius manifolds exists.
Furthermore, we define the notion of a tensor product diagram 
of Frobenius manifolds with factorizable flat identity
and prove the existence such a diagram and hence a tensor product
Frobenius manifold. These diagrams and manifolds are unique up
to equivalence. Finally, we
derive the special initial conditions for a tensor product of semi--simple
Frobenius manifolds in terms of the special initial conditions of the factors.
\end{abstract}

\maketitle

\section*{Introduction}

This paper is devoted to the study of Frobenius manifolds and their
tensor products.

\sm

The foundations of the theory of Frobenius manifolds were
laid down  by Dubrovin [D1].
Such manifolds play a central role
in the study of 
quantum cohomology and mirror symmetry (cf.\ [Gi, KM, M3]).
In the realm of mathematical physics 
they appear for instance as the canonical moduli spaces
of Topological Field Theories (TFTs).
Recently they have also emerged in the study of differential 
Gerstenhaber--Batalin--Vilkovisky algebras [BK].
Frobenius manifolds first arose in Saito's study of unfolding
singularities (cf.\ [D2, M3]) where they were called flat structures [S].
Further examples of Frobenius manifolds built on extended
affine Weyl groups were constructed in [DZh]. 
For an introduction to the subject the reader can
consult [D2, H, M1].
\sm

In a way one can regard Frobenius manifolds as a non--linear structure
on cohomology spaces. This non--linear structure is rigid and has weak
functorial properties. However, it admits a remarkable tensor product
operation.

\sm

In the formal setting this operation has been introduced in [KM].
In quantum cohomology it corresponds to the K\"unneth formula  [B, K].  

\sm

The main goal of this paper is to introduce and to study the
tensor product of analytic Frobenius manifolds in the non--formal
setting. This is done in two steps. First, we show that the formal tensor
product of two convergent potentials is convergent, so that we can define the 
tensor product of two pointed germs of analytic Frobenius manifolds.
The convergence proof is based upon the study of 
one--dimensional Frobenius manifolds, carried out in [KMZ].
Secondly, in the presence
of flat identities, we show that inside the convergence domain the so defined 
tensor products corresponding to different base--points are canonically
isomorphic. This observation is translated into the existence 
of a natural affine tensor product connection on the exterior product
of the tangent bundles over the Cartesian product of two
Frobenius manifolds. Using this connection 
we define the notion of 
tensor product diagrams for  Frobenius manifolds
with factorizable 
flat identities. To patch together the local pointed tensor products
we need the technical assumption that the flat identities of the factors
are factorizable which means that they can be split off as a Cartesian factor
${\Bbb C}$. This is the case in all important examples. 
In this situation, we prove the existence of such diagrams and
uniqueness up to equivalence. One of the pieces of data for these diagrams
is a tensor product 
Frobenius manifold for two Frobenius manifolds,
which contains a submanifold parameterizing all possible tensor products. The
size of the manifold itself depends on the convergence domain of the
tensor product potentials and cannot be controlled a priori. Therefore,
we regard two tensor product manifolds as equivalent
if they agree in an open neighborhood of this submanifold.
We also extend the situation to a slighly more general setting and
prove the analogous results.

\sm

In the examples stemming from the unfolding of singularities
the tensor product corresponds to the direct sum of singularities [M3].
For TFTs it provides the canonical
moduli space for the tensor product of two such theories.
The theorem of the existence of a tensor product then implies
that the moduli space obtained by tensoring all
possible natural perturbations of two given TFTs is included as
a subspace in the natural moduli space of the tensor theory.

\sm

Frobenius manifolds often carry the additional structures of an
Euler field and a flat identity
which are sometimes included in the definition [D2].
Our tensor product can also be extended into this category.
\sm

In the special case of semi--simplicity 
the structure of Frobenius manifolds becomes particularly
transparent.
Roughly speaking, the Frobenius structure is
determined by the Schlesinger special initial conditions [D2, M1, MM]
at a given tame base--point. In this setting we calculate the special
initial conditions of the tensor product.

\sm

Since we need to review the formalism of formal Frobenius manifolds
and their tensor product, this paper gives a complete analysis
of the tensor product in the theory of Frobenius manifolds in all of
its presently known facets.

\sm 

The paper is organized as follows:
We begin by recalling the necessary definitions and facts of the
theory of formal Frobenius manifolds and Frobenius manifolds, including
the tensor product in the formal setting in section 1. We also introduce
the notion of pointed germs of analytic 
Frobenius manifolds and give a one--to--one
correspondence with convergent formal Frobenius manifolds.
In section 2 we define the tensor product in the category of 
formal Frobenius manifolds with flat identity and Euler field and 
prove that the tensor product of two convergent potentials is again 
convergent yielding a tensor product for pointed germs of Frobenius manifolds.
Section 3 contains the definition of a global version of the tensor product
in terms of an affine tensor product connection on the exterior product of
the tangent bundles over the Cartesian product of 
two Frobenius manifolds. In the framework
of tensor product diagrams for  
Frobenius manifolds with factorizable flat identities,
we prove the existence of such a diagram  and show
uniqueness up to equivalence. By generalizing the
setting to general tensor product diagrams and introducing 
new natural restrictions, 
we are again able to show existence and uniqueness up to equivalence. 
The last section is an application of the previous results to semi--simple
Frobenius manifolds. In particular, we calculate the special
initial conditions at a tame semi--simple point
of the tensor product
Frobenius manifold in terms of those of the pre--images.

\sm

\subsection*{Acknowledgments}
I gratefully thank Yu.\ Manin for many invaluable discussions.
I would also like to thank B.\ Dubrovin for stimulating 
conversations on the subject and W.\ Kramer for 
his helpful remarks.
Finally, I thank the Max--Planck--Gesellschaft for financial support.

\section {\label{st1}Frobenius manifolds}

\sm
We begin by reviewing the necessary material from the theory of Frobenius 
manifolds:

\sm

\subsection{\label{ts1.1}Formal Frobenius manifolds}

\sm

We will follow the definition from [M1]. Let $k$ be a supercommutative
$\Bbb Q$--algebra, $H= \oplus_{a \in A} k \del_a$ a free (${\Bbb Z}_2$--graded)
$k$--module of finite rank, $g: H \otimes H \rightarrow k$ an
even symmetric pairing which is non--degenerate in the sense that
it induces an isomorphism $g': H \rightarrow H^t$ where $H^t$ is the dual 
module.

\sm
Denote by $K= k[[H^t]]$ the completed symmetric algebra of $H^t$.
This means that if $\sum_a x^a \del_a$ is a generic element of $H$, then
$K$ is the algebra of formal series $k[[x^a]]$. We will also regard
elements of $K$ as derivations on $H_K:= K\otimes_k H$ 
with $H$ acting via contractions.
We will call the elements of $H$ {\it flat}.

\sm

\subsubsection{\label{ts1.2}Definition}
The structure of {\em a formal Frobenius manifold on $(H,g)$} is
given by a potential $\Phi \in K$ defined up to quadratic terms
which satisfies the associativity of WDVV--equations:
\begin{equation}
\forall a,b,c,d: \sum_{ef}\Phi_{abe}g^{ef}\Phi_{fcd}=
(-1)^{\tilde a(\tilde b+\tilde c)} \sum_{ef}\Phi_{bce}g^{ef}\Phi_{fad}
\label{t1.1}
\end{equation}
where $\Phi_{abc} =\del_a\del_b\del_c\Phi$, $g^{ij}$ is the inverse metric
and $\tilde a := \widetilde {x^a}= \widetilde {\del_a}$ is 
the ${\Bbb Z}_2$--degree.

\sm
From the equations (\ref{t1.1}) it follows that the multiplication law
given by $\del_a \circ \del_b = \sum_c \Phi_{ab}^c \del_c$ 
turns $H_K = K \otimes_k H$ into a supercommutative $K$--algebra.

\bigskip

There are two other equivalent descriptions of formal Frobenius manifolds
using abstract correlation functions and $Comm_{\infty}$--algebras (cf.\ [M1]).

\sm

\sm

\subsubsection{\label{ts1.3}Definition}
{\em An abstract tree level system of correlation functions (ACFs)
on $(H,g)$} is a family of  $S_n$--symmetric even polynomials
\begin{equation}
Y_n: H^{\otimes n} \rightarrow k, \; n\geq 3
\label{t1.2}
\end{equation}
satisfying the Coherence axiom (\ref{t1.3}) below.

\sm
Set $\Delta = \sum \del_a g^{ab}\del_b$. Choose
any pairwise distinct $1 \leq i,j,k,l \geq n$ and  denote by
$ijSkl$ any 
partition $S = \{S_1,S_2\}$ of $\{1,\dots,n\}$ which 
separates $i,j$ and $k,l$, i.e.
$i,j \in S_1$ and  $k,l \in S_2$. The axiom now reads: 
\begin{itemize}
\item [] {\it Coherence:} For any
choice
of $i,j,k,l$
\begin{multline}
\sum_{ijSkl} \sum_{a,b} Y_{|S_1|+1} 
(\bigotimes_{r\in S_1} \gamma_r \otimes \del_a) g^{ab} 
Y_{|S_2|+1} (\del_b \otimes \bigotimes_{r\in S_2}
\gamma_r)\\
=  \sum_{ikTjl}  \sum_{a,b} Y_{|T_1|+1} 
(\bigotimes_{r\in T_1} \gamma_r \otimes \del_a) g^{ab} 
Y_{|T_2|+1}(\del_b \otimes \bigotimes_{r\in T_2} \gamma_r).
\label{t1.3}
\end{multline}
\end{itemize}

\sm

\subsubsection{\label{ts1.4}Correspondence between formal series and  
families of polynomials}

\sm

Given a formal series $\Phi \in K$, we can expand it up to terms of order
two as
\begin{equation}
\Phi = \sum_{n\geq 3}^{\infty} \frac{1}{n!} Y_n
\label{t1.4}
\end{equation}
where the $Y_n\in (H^t)^{\otimes n}$.  We will consider the
$Y_n$ as even symmetric maps $H^{\otimes n}\rightarrow k$.
One can check that the WDVV--equations (\ref{t1.1}) and the Coherence
axiom (\ref{t1.3}) are equivalent under this identification, see e.g.\ [M1].

\sm

\subsubsection {Remark}
Using  $Y_{n}$ one can define
multiplications  $\circ_n$ by dualizing with $g$ 
\begin{equation}
g(\circ_n(\gamma_1, \dots , \gamma_n),\gamma_{n+1}):= 
Y_{n+1}: H^{\otimes(n+1)} \ra k, \, Y_{n+1}(\gamma_1 \otimes \dots \otimes 
\gamma_n)
\label{acfmult}
\end{equation}
which define a so--called $Comm_{\infty}$--algebras.

\sm

\subsubsection{\label{ts1.5}Theorem (III.1.5 of [M1])}
{\it
The correspondence of \ref{ts1.4} establishes a bijection between the 
following structures on $(H,g)$.
\begin{itemize}
\item [(i)] Formal Frobenius manifolds.
\item [(ii)] Cyclic $Comm_{\infty}$--algebras.
\item [(iii)] Abstract correlation functions.
\end{itemize}
}
\sm

\subsubsection{\label{ts1.6}Definition}
An even element $e$ in $H_K$ is called
an {\em identity}, if it is an identity for the multiplication $\circ$.
It is called flat, if $e \in H$. In this case, we will denote
$e$ by $\del_0$ and include it as a basis element.

\sm

\subsubsection{\label{ts1.7}Euler Operator}

\sm

An even element $E\in K$ is called {\it conformal}, 
if ${\mathrm Lie}_E (g) = D g$ for 
some $D \in k$. Here, we take the Lie derivative of the 
tensor $g$ bilinearly extended to $K$
w.r.t.\ the derivation $E$. In other words:
\begin{equation}
\forall X,Y \in K: \quad {\mathrm Lie}_E (g) :=Eg(X,Y) -g([E,X],Y)- g(X,[E,Y])=
D\, g(X,Y).
\label{t1.5}
\end{equation}

It follows that $E$ is the sum of infinitesimal rotation, dilation and constant shift,
hence, we can write $E$ as:
\begin{equation}
E= \sum_{a,b \in A} d_{ab} x^a \del_b + \sum_{a \in A} r^a \del_a := E_1 + E_0,
\label{t1.6}
\end{equation}
for some $d_{ab} \in k$. 
Specializing $X= \del_a, Y=\del_b$ we can rewrite (\ref{t1.5})
\begin{equation}
\forall a,b:\quad \sum_c d_{ac}  g_{cb} +\sum_c d_{bc} g_{ac} = D g_{ab}.
\label{t1.7}
\end{equation}

In particular, we see that $[E,H] \subset H$ and that the operator
\begin{equation}
{\Cal V}: H \ra H: \quad {\Cal V}(X):= [X,E]- \frac{D}{2}X
\label{t1.8}
\end{equation}
is skew--symmetric.

\sm
A conformal operator $E$ is called {\it Euler}, if it additionally satisfies
${\mathrm Lie}_E(\circ) = d_0 \circ$ for some constant $d_0$.

\sm

\subsubsection{\label{ts1.8}Quasi--homogeneity}

\sm

The last condition is equivalent to the quasi--homogeneity
condition (Proposition 2.2.2. of [M1])
\begin{equation}
E\Phi = (d_0+D) \Phi + \text {a quadratic polynomial in flat
coordinates}.
\label{t1.9}
\end{equation}

\sm
\subsection{The tensor product of formal Frobenius manifolds}

\sm

The tensor product of formal Frobenius manifolds is naturally defined via
the respective Cohomological Field Theories (Cf.\ [KM, KMK, K])
which in terms of correlation
functions manifests itself in the appearance of operadic correlation
functions and the diagonal class $\Delta_{\mbn}\in A^{n-3}(\mbn\times \mbn)$.  

\sm

\subsubsection{\label{ts2.2}Trees and the cohomology of the spaces $\mbs$}

\sm

We will consider a tree $\tau$ as 
quadruple $(F_{\tau},V_{\tau},\del_{\tau},j_{\tau})$ of a (finite)
set of (of flags) $F_{\tau}$, a (finite) set (of vertices) $V_{\tau}$,
the boundary map $\del_{\tau}: F_{\tau} \ra V_{\tau}$,
and an involution $j_{\tau}  F_{\tau}\ra F_{\tau}, j_{\tau}^2=j_{\tau}$.

\sm
We call a a tree $S$--labeled or an $S$--tree if there is a fixed
isomorphism of the tails (one element orbits of $j$) and $S$.
We will only consider trees at least three tails.
A tree is called stable if the set of flags at each vertex is at least
of cardinality $3$: $\forall v \in V(\tau) |F_{\tau}(v)| \geq 3$.

\sm

If a tree $\tau$ is unstable we define the {\em stabilization} of
$\tau$ to be the tree obtained from $\tau$ by contracting one edge
at each unstable vertex. There are just three possible configurations
at a given unstable vertex
and it is easily seen that the result of the stabilization is indeed
stable and that the stabilization does not depend on the chosen edge.

\sm

\subsubsection {\label{ks1.2}Keel's presentation.}

\smallskip

As was shown in [Ke], 
the cohomology ring of $\mbs$ can be presented in terms of 
classes of boundary divisors as generators and quadratic 
relations as introduced by [Ke]. 
Thus we have a map  
\begin{equation}
[\; ]:\{\text{stable S--trees}\} \to H^*(M_{0S})
\label{ks1.3.1}
\end{equation}
The additive
structure of this ring and the respective
relations can then be naturally described in terms of 
stable trees (see [KM] and [KMK]).

\sm
\subsubsection {\label{ts2}Operadic Correlation Functions}

\sm

By identifying the index set $\overline n = \n$ of ACFs
or more generally any finite
set $S$
with a set of markings of a
$S$--tree, one can extend the notion of ACFs to operadic
correlation 
functions. These are maps from $H^{\otimes S}$ which also
depend on a choice of a stable $S$--tree $\tau$ (cf.\ [KM]).
\begin{equation}
Y(\tau) : H^{\otimes T_{\tau}} \rightarrow k
\label{t2.8}
\end{equation}
 
In fact under certain natural restrictions 
there exists a unique extension to trees
for any system of ACFs 
$\{Y_n\}$ (cf.\ Lemma 8.4.1 of [KM]).

\sm
\subsubsection{\label{ts2.5.2}Remark}
Given a set of ACFs $\{Y_n\}$ the correlation function of the
above cited Lemma for a
stable $n$--tree $\tau$ is given by the formula 
\begin{equation}
 Y(\tau)(\del_{a_1} \otimes \dots \otimes \del_{a_n})
= (\bigotimes_{v \in V_{\tau}} Y_{F_v}) 
(\del_{a_1} \otimes \dots \otimes \del_{a_n}\otimes
\Delta^{\otimes|E_{\tau}|}).
\label{t2.11}
\end{equation}

If fact due to the Coherence axiom (\ref{t1.3}) the operadic correlation
functions only depend on the class of the tree $[\tau]\in H^*(\mbn)$:
$Y(\tau) = Y([\tau])$.  

\sm

\subsubsection{\label{ts2.5.3}Remark}
To shorten the formulas, by abuse of notation, we will also
denote the following function from $H^{\otimes F_{\tau}}$ to $k$
by $Y(\tau)$:
\begin{equation}
\bigotimes_{v \in V_{\tau}} Y_{F_v} =: Y(\tau).
\label{t2.13}
\end{equation}

Which function is meant will be clear 
from the index set of the arguments.

\sm

\subsubsection{\label{ts3.1.1}The diagonal of $\mbn \times \mbn$}

\sm

Denote the class of the diagonal
in $H^{n-3}(\mbn \times \mbn)$ by $\Delta_{\mbn}$ and 
write it in terms of a tree basis:
\begin{equation}\Delta_{\mbn} = \sum_{[\sigma], [\tau]\in {\Cal B}_n}
[\sigma] g^{[\sigma] [\tau]}\otimes {[\tau]}
\label{t3.1}
\end{equation}
where 
$g_{[\sigma] [\tau]} = \int_{\mbn}[\sigma] \cup [\tau]$
and ${\Cal B}_n$ is some basis of $H^*(\mbn)$.
Notice that 
\begin{equation}
g^{[\sigma] [\tau]} =0
\; \text { unless } \; 
|E_{\sigma}|+|E_{\tau}|=n-3.
\label{t3.2}
\end{equation}

For explicit computations one can 
use the basis ${\Cal B}_n$ presented in [K].

\sm
\subsubsection {Tensor product for ACFs}
The tensor product of two systems of ACFs 
$(H^{(1)}$, $\Delta^{(1)}$, $\{Y^{(1)}_n\})$ and
$(H^{(2)}$, $\Delta^{(2)}$, $\{Y^{(2)}_n\})$ is the system of ACFs 
$(H^{(1)}\otimes H^{(2)}, \Delta^{(1)} \otimes \Delta^{(2)}, \{Y_n\})$ 
defined by
\begin{multline}
Y_n ((\gamma^{(1)}_1 \otimes\gamma^{(2)}_1) \otimes \dots
\otimes (\gamma^{(1)}_n \otimes \gamma^{(2)}_n)) :=\\
\epsilon (\gamma^{(1)},\gamma^{(2)})(Y^{(1)} \otimes Y^{(2)})(\Delta_{\mbn})
((\gamma^{(1)}_1 \otimes\gamma^{(2)}_1) \otimes \dots
\otimes (\gamma^{(1)}_n \otimes \gamma^{(2)}_n)).
\label{t3.4}
\end{multline}
where for each summand
\begin{multline}
(Y^{(1)} \otimes Y^{(2)})([\tau] \otimes [\sigma])
((\gamma^{(1)}_1 \otimes\gamma^{(2)}_1) \otimes \dots
\otimes (\gamma^{(1)}_n \otimes \gamma^{(2)}_n)) =\\ 
\epsilon (\gamma^{(1)},\gamma^{(2)})Y^{(1)}
({\tau})(\gamma^{(1)}_1 \otimes \dots
\otimes \gamma^{(1)}_n)
Y^{(2)}({\sigma})(\gamma^{(2)}_1 \otimes \dots \otimes \gamma^{(2)}_n).
\label{t3.3}
\end{multline}
\sm

\subsubsection{\label{ts3.1.5}Definition}
Given two formal Frobenius manifolds $(H^{(1)}$, $g^{(1)}$, $\Phi^{(1)})$ 
and $(H^{(2)}$, $g^{(2)}$, $\Phi^{(2)})$,
let $\{Y^{(1)}_n\}$ and $\{Y^{(2)}_n\}$ be the corresponding
ACFs. 
{\em The tensor product $(H,g,\Phi)$ of $(H^{(1)},g^{(1)},\Phi^{(1)})$ 
and $(H^{(2)},g^{(2)},\Phi^{(2)})$} is 
defined to be $(H^{(1)}\otimes H^{(2)}, g^{(1)}\otimes g^{(2)},\Phi)$ 
where the potential $\Phi$ 
is given by:
\begin{equation}
\Phi(\gamma) = \sum_{n\geq3} \frac{1}{n!}
(Y^{(1)}\otimes Y^{(2)})(\Delta_{\mbn})(\gamma^{\otimes n}).
\label{t3.5}
\end{equation}

\sm

As in \ref{ts1.4}, to make sense of (\ref{t3.5}) one should expand 
$\gamma = \sum x^{a',a''} \del_{a'a''}$ in 
terms of the tensor product 
basis $(\del_{a'a''}:= \del^{(1)}_{a'}\otimes \del^{(2)}_{a''})$
of the two basis  $\{\del^{(1)}_{a'}\}$ and $\{\del^{(2)}_{a''} \}$ 
and the dual coordinates $x^{a'a''}$ for this basis.

\bigskip

Inserting the explicit basis ${\Cal B}_n$ with its intersection form
given in [K] allows to make (\ref{t3.5}) explicit (cf.\ [K]). W 
In the case of quantum cohomology this provides the explicit 
K\"unneth formula. Here the tensor product
of the potentials $\Phi^V$ and $\Phi^W$
belonging to some smooth projective varieties $V$ and $W$
is the Gromov--Witten potential of $\Phi^{V\times W}$ (cf.\ [KM, KMK, K, B]).

\subsection {Frobenius manifolds}

\sm

\subsubsection{\label{ts4.1.1}Definition}

\sm
A {\em Frobenius manifold} $\bold M$
is a quadruple $(M,{\Cal T}_M^f,g,\Phi)$ of a 
(super)\-mani\-fold $M$, an affine flat structure ${\Cal T}_M^f$,
a compatible
metric $g$ and a potential function whose tensor of 
third derivatives defines an associative commutative
multiplication $\circ$
on each fiber of ${\Cal T}_M$.

\sm
For the notion of supermanifolds and supergeometry in general
we refer to the book [M2]. Below we will consider only
manifolds in the analytic category.
\sm

\subsubsection{\label{ts4.1.2}Definition} 
A {\it pointed Frobenius manifold} is a pair $({\bold M},m_0)$
of a Frobenius
manifold $\bold M$ and a point $m_0 \in M$ called the base--point. 

\sm
When considering {\it flat coordinates} in a neighborhood of 
the base--point $m_0$ of a pointed Frobenius manifold, we
require  
that the coordinates of $m_0$ are all zero. In other words,
the base--point
corresponds to a choice of a zero in flat coordinates.

\sm

\subsubsection{\label{ts4.1.3}Euler field and Identity}

\sm

Just as in the formal case, a
Frobenius manifold may carry two additional structures; an 
Euler field and an identity. They are defined analogously.

\sm
\subsubsection{Definition.} 
An even vector field $E$ 
on a Frobenius manifold with a flat metric $g$ is called
{\em conformal} of conformal weight $D$, for some constant $D$, 
if it satisfies ${\mathrm Lie}_E(g) = Dg$.
A conformal field $E$ is called {\it Euler}, if it additionally
satisfies 
${\mathrm Lie}_E(\circ) = d_0 \circ$ for some constant $d_0$.

\sm

\subsubsection{\label{ts4.1.4}From germs of pointed Frobenius manifolds to
convergent formal Frobenius manifolds}

\sm

Regarding a germ of a pointed Frobenius manifold $({\bf M},m_0)$
over a field $k$ of characteristic zero, choose a flat basis of vector fields
$(\del_a)$ and set $H= \oplus_a k \del_a$ and
keep the metric $g$. Choose corresponding unique
local flat coordinates $x^a$ s.t.
$\forall a: x^a(m_0)=0$ as we demanded in \ref{ts4.1.2}.
A structure of a formal Frobenius manifold on $(H,g)$ is then 
given 
by the expansion of the potential into a power series in
local flat coordinates $(x^a)$ at $m_0$. Up to 
quadratic terms we obtain:
\begin{equation}
\Phi_{m_0}({\bold x}) = \sum_{n \geq 3} \frac{1}{n!} 
\sum_{a_1, \dots ,a_n \in \n} x^{a_n} \cdots x^{a_1}
Y_n^{m_0}(\del_{a_1} \otimes \dots \otimes \del_{a_n})
\label{t4.1}
\end{equation}
where the functions $Y_n$ are defined via
\begin{equation}
Y_n^{m_0} (\del_{a_1} \otimes \dots \otimes \del_{a_n}) :=
\del_{a_1}  \cdots \del_{a_n} \Phi\vert_{m_0} 
=(\del_{a_1} \otimes \dots \otimes \del_{a_n}) \Phi_{m_0}({\bold x})\vert_0 
\label{t4.2}
\end{equation}
$\Phi_{m_0}$ obviously obeys the WDVV--equations.

\sm

Furthermore, in the presence of an Euler field
or a flat identity  writing $E$ and $e=\del_0$ in flat coordinates
defines the same structures in the formal situation.

\sm

We stress again that we are dealing with pointed
Frobenius manifolds. Due to this a zero in flat
coordinates has been fixed and $\{Y_n^{m_0}\}, E$ and $e$ are
uniquely defined.

\sm
On the other hand, the functions
in (\ref{t4.2}) and $E$ are dependent on the choice of the base--point.
Choosing a different base--point $\widehat {m}_0$ with
$x$--coordinates
$x^a(\widehat {m}_0)= x_0^a$ in the domain of
convergence of $\Phi_{m_0}$ yields the new standard flat coordinates
$\widehat{x}^a = x^a - x_0^a$. The
corresponding 
functions $Y$ transform via:
\begin{multline}
\widehat{Y}^{\widehat {m}_0}_n (\del_{a_1} \otimes \dots \otimes \del_{a_n}) :=
 \del_{a_1}  \dots \del_{a_n} \Phi\vert_{\widehat {m}_0}
=  \del_{a_1}  \dots \del_{a_n} \Phi_{m_0}\vert_{x^a_0}\\
= \sum_{N \geq 0}\frac {1}{N!} \sum_{(b_1,\dots,b_N):b_i \in A}
x^{b_N}_0 \cdots x^{b_1}_0
Y^{m_0}_{n+N}(\del_{b_1} \otimes \dots \otimes \del_{b_N} \otimes
\del_{a_1} \otimes \dots \otimes \del_{a_n})\\
= \sum_{N \geq 0}\frac {1}{N!} \sum_{(b_1,\dots,b_N):b_i \in A}
\eps(b|a)\,
x^{b_N}_0 \cdots x^{b_1}_0
Y^{m_0}_{n+N}(\del_{a_1} \otimes \dots \otimes \del_{a_n} \otimes
\del_{b_1} \otimes \dots \otimes \del_{b_N}) 
\label{t4.3}
\end{multline}
where $\eps(b|a)$ is a shorthand notation
for $\eps(\del_{b_1} \cdots \del_{b_N}|\del_{a_1}\dots \del_{a_n })$
which we define as the superalgebra sign 
acquired by permuting $\del_{b_1},\dots , \del_{b_N}$ 
past the $\del_{a_1},\dots , \del_{a_n}$:
\begin{equation}
\del_{b_1} \cdots \del_{b_N}\del_{a_1}\cdots \del_{a_n}
= \eps(b|a)\,\del_{a_1}\cdots\del_{a_n}\del_{b_1} \cdots \del_{b_N}.
\end{equation}

\subsubsection{Notation}
We denote the convergent formal Frobenius structure obtained from a pointed
Frobenius manifold $((M,T_M^f,g_M,\Phi^M),p)$ 
with a choice of a basis $(\del_a)$ of $T_M^f$
by
\begin{equation}
(\bigoplus k \; \del_a,g_M,\Phi^M_p)
\end{equation}

\sm

\subsubsection{\label{ts4.1.5}From convergent formal Frobenius manifolds 
to germs of poin\-ted Frobenius manifolds}

\sm

Starting with any formal Frobenius manifold $(H,g)$ 
with a potential $\Phi$, we can produce a germ of a manifold with
a flat structure by identifying the $x^a$ as coordinate functions
around some point $m_0$, choosing $H$ as the space of
flat fields and considering $g$ as the metric. To get a Frobenius
manifold, however, we need that the formal potential $\Phi$ has
some nonempty domain of convergence. If
\begin{equation}
\Phi(\gamma) = \sum_{n\geq3} \frac {1}{n!} Y_n(\gamma^{\otimes n})
\label{t4.4}
\end{equation}
with $\gamma = \sum x^a \Delta_a$ is convergent, we can pass
to a germ of a pointed Frobenius manifold. If necessary,
we can, in this situation, even move the base--point as indicated above.

\section{The Tensor product for Euler fields, flat identities and
germs of Frobenius manifolds}

\subsection{\label{ts3.2}The tensor product for Euler fields
and flat identities} 

\med

In this section, we extend the operation of forming the tensor product
to the additional structures of an Euler field and an identity.
In order to achieve this, we first rewrite the quasi--homogeneity
condition and the defining relation for an identity in terms of
operadic correlation functions. To this end we introduce
the morphisms $\pi^*,\pi_*$ on trees.

\subsubsection{\label{ts2.3}Forgetful morphisms and trees.}

\sm

The flat and proper morphisms 
$\pi_s: \mbs \rightarrow \Mb_{0,S\setminus \{s\}}$ which
forget the
point marked by $s$ and stabilize 
if necessary induce the maps  
$\pi_*$ and $\pi^*$
on the Chow rings where we omitted the subscript $s$ which we
will always do, if there is no risk of confusion.

\sm
We will now define the maps $\pi_*,\pi^*$ on trees
corresponding
under $[\;] $ to the respective maps in the Chow rings of $\mbs$.

\sm
Define $\pi_*$ via
\begin{align}
\hskip -15pt \pi_{s*}(\tau) =& \cases \text {forget the tail number $s$ and
stabilize, if the
stabilization is necessary}\hskip -1cm\\0 \; \text {otherwise} 
\endcases
\label {t2.5}\\
\intertext {For any $S$--tree $\tau$ and any $s \notin S$ set}
\tau^{s}_v = &\quad
\text {the ($S \cup \{s\}$)--tree obtained from $\tau$ by
adding an 
additional tail}\nn\\
&\quad \text{marked by $s$ at the vertex $v$}.
\end{align}
 
\sm

For the notion of stabilization of a tree cf.\ \ref{ts2.2}.

\sm
 
Now we define:
\begin{equation}
\pi^{s*}(\tau) = \sum_{v \in V_{\tau}} \tau_v^{s}.
\label {t2.6}
\end{equation}

\sm

Taking the definition of $[\; ]$ from [KM] 
it is a straightforward calculation using e.g.\ [Ke] to check
that indeed 
$[\pi^*(\tau)] = \pi^*([\tau])$ and 
$[\pi_*(\tau)] = \pi_*([\tau])$.

\med

\subsubsection{\label{ts2.6}Quasi--homogeneity condition in terms 
of correlation functions.}

\med 

\subsubsection{\label{ts2.6.1}Lemma.}
{\it
In terms of the abstract correlation functions $Y_n$ the
quasi--homogeneity condition (\ref{t1.9})
is given by
\begin{multline}
\sum_{a \in A} 
(\sum_{i=1}^n d_{a_i a} Y_n(\del_{a_1} \otimes \dots \otimes \widehat {\del_{a_i}}
\otimes \dots\otimes
\del_{a_n}\otimes \del_a) +
r^a Y_{n+1} (\del_{a_1} \otimes \dots \otimes
\del_{a_n} \otimes \del_a))\\
= (d_0+D) Y_n(\del_{a_1} \otimes \dots \otimes \del_{a_n}).
\label{t2.14}
\end{multline}
}

\sm

{\bf Proof.} Applying the vector field $E$ in the form (\ref{t1.6}) to
(\ref{t1.4}) and making a coefficient check yields (\ref{t2.14}).

\med

\subsubsection{\label{ts2.6.2}Lemma.}
{\it
The correlation functions (\ref{t2.13}) obey the following 
relation. For a given $n$--tree $\tau$: 
}
\begin{multline}
\sum_{a \in A}
(\sum_{f  \in F_{\tau}} d_{fa}  Y(\tau) 
((\bigotimes_{f'  \in F_{\tau}\setminus\{f\}}\del_{f'}) \otimes \del_a) 
+  r^a Y(\pi^*(\tau))((\bigotimes_{f  \in F_{\tau}}
\del_f) \otimes \del_a))\\
=|V_{\tau}| (D+d_0) Y(\tau) (\bigotimes_{f  \in F_{\tau}} \del_f).
\label{t2.15}
\end{multline}

\sm

{\bf Proof.} Recall that by definition 
$Y(\pi^*(\tau)) = \sum_{v\in V_{\tau}} Y(\tau_v^{n+1})$. 
By applying (\ref{t2.14}) at every vertex $v$ of $\tau$, we obtain
\begin{equation}
\begin{split}
\sum_{f  \in F_{\tau}}& \sum_a d_{fa}\, Y(\tau) 
((\bigotimes_{f'  \in F_{\tau}\setminus\{f\}}\del_f) \otimes \del_a) 
\hskip 6cm\\
&=\sum_{v \in V_{\tau}} [\sum_{f\in F_{\tau}(v)} d_{fa} 
\bigotimes_{v' \in V_{\tau}}
(Y_{|F_{\tau}(v')|})((\bigotimes_{f'  \in F_{\tau}\setminus\{f\}}\del_{f'}) 
\otimes \del_a) ]
\\
&= \sum_{v \in V_{\tau}}[ 
(D+d_0)Y(\tau) (\bigotimes_{f  \in F_{\tau}} \del_f)  
- \sum_{a\in A} r^a Y(\tau_{v}^{n+1})((\bigotimes_{f  \in
F_{\tau}} \del_f) \otimes \del_a)]\\
&= |V_{\tau}| (D+d_0) Y(\tau) (\bigotimes_{f  \in F_{\tau}}
\del_f) -  \sum_{a \in A} r^a
Y(\pi^*(\tau))((\bigotimes_{f  \in F_{\tau}} \del_f) \otimes
\del_a).
\end{split}
\end{equation}

\med

\subsubsection{\label{ts2.6.3}Proposition}
{\it
For the operadic correlation functions $\{Y(\tau)\}$ the
quasi--homogeneity
condition is equivalent to
}
\begin{multline}
\sum_{i=1}^n \sum_{a\in A} d_{a_i a} Y(\tau) (\del_{a_1}
\otimes \dots\otimes \widehat{\del_{a_i}} \otimes\dots \otimes \del_{a_n}
\otimes \del_a) 
-|E_{\tau}|\,d_0 Y(\tau)(\del_{a_1} \otimes \dots \otimes \del_{a_n})\\
+  \sum_{a \in A} r^a Y(\pi^*(\tau)) (\del_{a_1} \otimes \dots
\otimes \del_{a_n} \otimes \del_a) 
= (d_0+ D) Y(\tau)(\del_{a_1} \otimes \dots \otimes
\del_{a_n}).
\label{t2.16}
\end{multline}

\sm

{\bf Proof.}
Writing out  the Casimir elements $\Delta = \sum \del_p g^{pq}\del_q$ 
and applying
Lemma \ref{ts2.6.2} yields:
\begin{multline}
\hskip -15pt{\shoveleft \sum_{i=1}^n \sum_{a\in A} d_{a_ia}  Y(\tau)(\del_{a_1}
\otimes \dots \otimes \widehat {\del_{a_i}} \otimes \dots
\otimes \del_{a_n} \otimes \del_a)+ |E_{\tau}|\, D \, Y(\tau)(\del_{a_1}
\otimes \dots 
\otimes \del_{a_n})} \\
\shoveleft{= \sum_{i=1}^n \sum_{a\in A} d_{a_ia} 
Y(\tau)(\del_{a_1}
\otimes \dots \otimes \widehat {\del_{a_i}} \otimes \dots
\otimes \del_{a_n} \otimes \del_a)}\\
+ |E_{\tau}|\, D 
\sum_{{(p_1,\dots p_{|E_{\tau}|}) p_i\in A}\atop
{(q_1,\dots q_{|E_{\tau}|}) q_i\in A}} 
(\bigotimes_{v \in V_{\tau}} Y_{F_v}) 
(\del_{a_1} \otimes \dots \otimes \del_{a_n}\otimes
\bigotimes_{j=1}^{|E_{\tau}|} (\del_{p_j} g^{p_j q_j}\otimes
\del_{q_i}))\\
\shoveleft {\overset {(*)}{=} 
\sum_{(p_1,\dots p_{|E_{\tau}|}) p_i\in A\atop
(q_1,\dots q_{|E_{\tau}|}) q_i\in A} [
\sum_{i=1}^n \sum_{a\in A} d_{a_ia}  Y(\tau)(\del_{a_1}
\otimes \dots \otimes \widehat {\del_{a_i}} \otimes \dots
\otimes \del_{a_n} \otimes \del_a \otimes
\bigotimes_{j=1}^{|E_{\tau}|}(\del_{p_j} g^{p_j q_j}\otimes
\del_{q_i}))}\\
 + \sum_{i=1}^{|E_{\tau}|} \sum_{a\in A}  d_{p_i a}(\bigotimes_{v \in V_{\tau}} Y_{F_v}) 
(\del_{a_1} \otimes \dots \otimes \del_{a_n}\otimes
\bigotimes_{j=1, j\neq i}^{|E_{\tau}|} (\del_{p_j} g^{p_j q_j}\otimes
\del_{q_i})\otimes  \widehat{\del_{p_i}}g^{p_ia}\del_{q_i}\otimes \del_a)\\
+ \sum_{i=1}^{|E_{\tau}|}
\sum_{a\in A}  d_{p_j a}(\bigotimes_{v \in V_{\tau}} Y_{F_v}) 
(\del_{a_1} \otimes \dots \otimes \del_{a_n}\otimes
\bigotimes_{i=1,i\neq j}^{|E_{\tau}|} (\del_{p_i} g^{p_i q_i}\otimes
\del_{q_i})\otimes  \del_{p_i}g^{p_ia}\widehat{\del_{q_i}}\otimes \del_a)] \\
\shoveleft{=(|E|_{\tau}+1)(D+d_0)Y(\tau)(\del_{a_1} \otimes \dots \otimes
\del_{a_n})}\\
-\sum_{a \in A} r^a  Y(\pi^*(\tau))
(\del_{a_1} \otimes \dots \otimes \del_{a_n}\otimes \del_a).
\label{t2.17}
\end{multline}

The equality $(*)$ holds due to (\ref{t1.7}). Rewriting (\ref{t2.17}),
we obtain (\ref{t2.16}). Vice versa postulating (\ref{t2.16}), we see that it
reduces to (\ref{t2.14}) for the one--vertex tree $(\rho_n)$.

\med

\subsubsection{\label{ts2.7}The identity in terms of correlation functions}

\sm

As previously remarked, we will assume that
the identity is a flat vector field
$e=\del_0$. As the semi--simplicity of $E$ this restriction
is satisfied in the case of quantum cohomology.

\med

\subsubsection{\label{ts2.7.1}Remark}

\sm

From Corollary 2.1.1 of [M1], we have that
\begin{equation}
Y_3(\del_a,\del_b,\del_0) = g_{ab} \text { and } 
Y_n (\del_{a_1} \otimes \dots \otimes \del_{a_{n-1}}\otimes
\del_0)
= 0
\quad \forall n >3 
\label{t2.19}
\end{equation}
are equivalent to the fact that $\del_0$ is a flat identity.

\med

\noindent In terms of operadic ACFs one obtains:

\sm

\subsubsection{\label{ts2.7.2}Proposition}
{\it 
For a flat identity $e= \del_0$ and for any stable $n$--tree
$\tau$ with $n > 3$}
\begin{equation}
Y(\tau)(\del_{a_1} \otimes \dots \otimes \del_{a_{n-1}}\otimes
\del_0) =
Y(\pi_*(\tau))(\del_{a_1} \otimes \dots \otimes \del_{a_{n-1}}).
\label{t2.20}
\end{equation}

\sm
{\bf Proof.} From (\ref{t2.19}) we know that 
$Y(\tau)(\del_{a_1} \otimes \dots \otimes \del_{a_{n-1}}\otimes \del_0) =0$,
if the valence of the vertex $v_0$ with the tail marked with $n$
is greater 
than three or, in other words, if the vertex remains stable after
forgetting the
tail $n$. Assume now that the vertex has valence three. Noticing
that
for a flat identity $Y_3(\del_a,\del_b,\del_0)= g_{ab}$ the
result follows
by direct calculation. There are two cases: either $v_0$ has two
tails marked
$n$ and $i$ for some $i$ and is joined to one other vertex $v'$
by the 
edge $e$ or $v_0$ just has one tail and is joined to two
other vertices
by the edges $e_1$ and $e_2$. In the first case we get
\begin{align*}
Y(\tau)&(\del_{a_1} \otimes \dots \otimes \del_{a_{n-1}} \otimes
\del_0 \otimes
\Delta^{\otimes |E_{\tau}|}) \hskip 5cm\\
&=(\bigotimes_{v \in V_{\tau} \setminus \{v_0\} }
(\bigotimes_{f \in F_{\tau}(v)} Y_{F_{\tau}(v)})\otimes
Y_{F_{\tau}(v_0)})\\
&\quad (\del_{a_1} \otimes \dots \otimes \widehat{\del_{a_i}} \otimes
\dots \otimes 
\del_{a_{n-1}} \otimes \Delta^{\otimes |E_{\tau}|-1}\otimes
\Delta_e  \otimes \del_{a_i}
\otimes \del_0)\\
&= \sum_{pq} (\bigotimes_{v \in V_{\tau}\setminus \{v_0\}}
(\bigotimes_{f \in F_{\tau}(v)}Y_{F_{\tau}(v)}))\\
&\quad (\del_{a_1} \otimes \dots \otimes \widehat{\del_{a_i}} \otimes
\dots \otimes \del_{a_{n-1}}
\otimes \Delta^{\otimes |E_{\tau}|-1}\otimes
\del_p g^{pq} g_{qa_i})\\
&= (\bigotimes_{v \in V_{\tau}\setminus \{v_0\}}
(\bigotimes_{f \in F_{\tau}(v)}Y_{F_{\tau}(v)}))
(\del_{a_1}\otimes \dots \otimes \del_{a_{n-1}} \otimes 
\Delta^{\otimes |E_{\tau}|-1})\\
&= Y(\pi_*(\tau)) (\del_{a_1} \otimes \dots \otimes
\del_{a_{n-1}})\\
\intertext {likewise in the second case}
Y(\tau)&(\del_{a_1} \otimes
\dots \otimes \del_{a_{n-1}} \otimes \del_0 \otimes
\Delta^{\otimes |E_{\tau}|}) \hskip 5cm\\
&=(\bigotimes_{v \in V_{\tau}\setminus \{v_0\}}(\bigotimes_{f \in
F_{\tau}(v)}
Y_{F_{\tau}(v)})\otimes Y_{F_{\tau}(v_0)})\\
& \quad (\del_{a_1} \otimes \dots \otimes \del_{a_{n-1}} 
 \otimes \Delta^{\otimes |E_{\tau}|-2}\otimes
\Delta_{e_1} \otimes \Delta_{e_2} \otimes \del_0)\\
&= \sum_{pq,rs} (\bigotimes_{v \in V_{\tau}\setminus \{v_0\}}
(\bigotimes_{f \in F_{\tau}(v)}Y_{F_{\tau}(v)}))\\
& \quad (\del_{a_1} \otimes \dots \otimes \del_{a_{n-1}} \otimes 
\Delta^{\otimes |E_{\tau}|-2}\otimes
\del_p g^{pq} g_{qr} g^{rs} \otimes \del_s)\\
&= (\bigotimes_{v \in V_{\tau}\setminus \{v_0\}}
(\bigotimes_{f \in F_{\tau}(v)}Y_{F_{\tau}(v)}))
(\del_{a_1} \otimes \dots \otimes \del_{a_{n-1}} \otimes 
\Delta^{\otimes |E_{\tau}|-1})\\
&= Y(\pi_*(\tau)) (\del_{a_1} \otimes \dots \otimes
\del_{a_{n-1}}).
\end{align*}

\med

\subsubsection{\label{ts2.7.3}Remark}

\sm 

In the setting of operads and higher order multiplications 
([G, GK]), the formulas (\ref{t2.19}) 
for a flat identity $e=\del_0$
correspond to the statements that $e$ is an identity  for
$\circ_2$ and acts as a
zero for all higher multiplications $\circ_n, n\geq 3$. The
contents of Proposition \ref{ts2.7.2} is the extension of these properties to any concatenation
of these multiplications. 

\med

After these preparations, we come to the main result of this section:
\med

\subsubsection{\label{ts3.2.1}Theorem}
{\it
Given two formal Frobenius manifolds $(H^{(1)},g^{(1)},\Phi^{(1)})$ and \hfill \linebreak
$(H^{(2)},g^{(2)},\Phi^{(2)})$
with Euler fields
\begin{align}
E^{(1)}&= \sum_{a'b' \in A^{(1)}} d^{(1)}_{a'b'} x^{(1) a'} \del^{(1)}_{b'} + \sum_{a' \in A^{(1)}}
r^{(1) a'} \del^{(1)}_{a'} \hskip 1cm
\text {of weight $D^{(1)}$ and}\\
E^{(2)}&= \sum_{a''b'' \in A^{(2)}} d^{(2)}_{a''b''} x^{(2) a''} \del^{(2)}_{b''} +
\sum_{a'' \in A^{(2)}} r^{(2) a''} \del^{(2)}_{a''} \hskip 1cm
\text {of weight $D^{(2)}$}
\end{align}
and with flat identities 
$e^{(1)},e^{(2)} $ of the same weight
$d^{(1)}_0= d^{(2)}_0=d$,
then
\begin{equation}
e = e^{(1)} \otimes e^{(2)} =\del^{(1)}_0 \otimes \del^{(2)}_0 = \del_{00}
\label{t3.6}
\end{equation}
and
\begin{multline}
E = \sum_{(b', b'') \in A^{(1)} \times A^{(2)}} [
\sum_{a' \in A^{(1)}} (d^{(1)}_{a'b'} x^{a'b''}) 
+ \sum_{a'' \in A^{(2)}} 
(d^{(2)}_{a''b''} x^{b'a''}) 
-  d x^{b'b''}]\; \del_{b'b''}\\
+ \sum_{a' \in A^{(1)}} r^{(1) a'} \del_{a'0}
+ \sum_{a'' \in A^{(2)}} r^{(2) a''} \del_{0a''}
\label{t3.7}
\end{multline}
define a flat identity of weight $d$ and an Euler field of weight
$\smash{D^{(1)}+D^{(2)}-2d}$
on the tensor product $(H,g,\Phi)$ of $(H^{(1)},g^{(1)},\Phi^{(1)})$ and $(H^{(2)},g^{(2)},\Phi^{(2)})$.
}

\med
Before we can prove the above theorem, we need one more Lemma about the
properties of the diagonal class $\Delta_{\mbn}$.

\med

\subsubsection{\label{ts3.2.2}Lemma}
{\it  
\begin{equation}
(id,\pi_*) (\Delta_{\mbn}) = (\pi^*,id)
(\Delta_{\overline{M}_{0n-1}})
\label{t3.8}
\end{equation}
and 
\begin{equation}
(\pi_*,\pi_*)(\Delta_{\mbn}) =0.
\label{t3.9}
\end{equation}
}
\med

{\bf Proof.}  Consider two any strata classes
$D_{\tau}\in A^*(\mbn), D_{\sigma} \in A^*(\mbnn)$.
Using the projection formula twice, we obtain
\begin{align*}
\int_{\mbn} D_{\tau} \cup \pi^*(D_{\sigma}) &= 
\int_{\mbnn} \pi_*(D_{\tau}) \cup D_{\sigma}\\
\Leftrightarrow \int_{\mbn \times \mbn}  (D_{\tau} \boxtimes
\pi^*(D_{\sigma}))
\cup \Delta_{\mbn}
&= \int_{\mbnn \times \mbnn} (\pi_*(D_{\tau})\boxtimes
D_{\sigma}) 
\cup \Delta_{\mbnn}\\
\Leftrightarrow \int_{\mbn \times \mbnn}  \hskip -7pt (D_{\tau} \boxtimes
D_{\sigma}) 
\cup (id,\pi_*)\Delta_{\mbn}
&=\int_{\mbn \times \mbnn} \hskip -7pt (D_{\tau} \boxtimes D_{\sigma}) 
\cup (\pi^*,id)\Delta_{\mbnn}.
\end{align*}
Since the intersection pairing is non--degenerate and the 
classes $D_{\tau} \boxtimes D_{\sigma}$ generate 
$A^*(\mbnn \times \mbn)$, the formula (\ref{t3.8})
follows.
Using the same type of argument for 
\begin{equation*}\begin{split}
\shave{\int_{\mbnn \times \mbnn}}&  (D_{\tau} \boxtimes
D_{\sigma}) 
\cup (\pi_*,\pi_*)\Delta_{\mbn} \\
&= \shave{\int_{\mbnn \times \mbnn}}  (\pi^*,\pi^*)(D_{\tau}
\boxtimes D_{\sigma}) 
\cup \Delta_{\mbn}\\
&=\int_{\mbn} \pi^*(D_{\tau}) \cup \pi^*(D_{\sigma}) = 
\int_{\mbn} \pi^*(D_{\tau} \cup D_{\sigma})=0
\end{split}
\end{equation*}
where the last zero is due to dimensional reasons, we obtain the
second claim
(\ref{t3.9}).

\med

{\bf Proof of the Theorem.}

\sm

As in \ref{ts3.1.5}, we choose the coordinates $x^{a'a''}$ corresponding
to the basis $\del_{a'}\otimes \del_{a''}$. 
The metric for the tensor product is
given by 
\begin{equation}
g_{a'b',a''b''}:= g (\del_{a'} \otimes \del_{a''},\del_{b'}
\otimes \del_{b''})=
g^{(1)}(\del_{a'} , \del_{b'}) g^{(2)}(\del_{a''} , \del_{b''})=
g^{(1)}_{a',b'}g^{(2)}_{a'',b''}.
\label{t3.10}
\end{equation}

\pagebreak

{\it Euler field.}

\med

First we check that $E$ is conformal of weight $D^{(1)}+D^{(2)}-2d$. On
the basis of
flat vector fields we calculate:
\begin{align}
&g([\del_{a'a''},E],\del_{b'b''})
+
g(\del_{a'a''},[\del_{b'b''},E])
\nn\\
&= \sum_{c'} d^{(1)}_{a'c'}g^{(1)}_{c' b'}g^{(2)}_{a''b''} 
+ \sum_{c''} d^{(2)}_{a''c''} g^{(1)}_{c'b'}g^{(2)}_{c''b''}
+ \sum_{c'}  d^{(1)}_{b'c'}g^{(1)}_{a'c'}g^{(2)}_{a''b''}
+ \sum_{c''} d^{(2)}_{b''c''}g^{(1)}_{a'b'}g^{(2)}_{a''c''}\nn\\
&\quad - 2 d g^{(1)}_{a'b'}g^{(2)}_{a''b''}\nn\\
&= (D^{(1)}+D^{(2)}-2d) g_{a'a'',b'b''}.
\label{t3.11}
\end{align}

We will prove the fact that $E$ is
indeed an Euler field by verifying the
quasi--homogeneity condition (\ref{t1.9}).

Set $D= D^{(1)}+D^{(2)}-2d$ and $\gamma = \sum x^{a'a''} \del^{(1)}_{a'}\otimes
\del^{(2)}_{a''}$:
\begin{align}
\hskip -1cm E_1&\Phi(\gamma)
= E_1 \sum_{n\geq 3} \frac{1}{n!} Y_n(\gamma^{\otimes n}) = 
E_1 \sum_{n\geq 3} \frac{1}{n!} 
(Y^{(1)}\otimes Y^{(2)})(\Delta_{\mbn})(\gamma^{\otimes n}) \hskip 2.5cm\nn\\
=& \sum_{n\geq 3} \frac{1}{n!}
 x^{a_n'a_n''}\cdots x^{a_1'a_1''} \Big(\sum_{a'\in A^{(1)}}\sum_{i=1}^n d^{(1)}_{a'_ia'} (Y^{(1)}\otimes Y^{(2)})(\Delta_{\mbn})\nn\\
&\qquad((\del^{(1)}_{a_1'}\otimes \del^{(2)}_{a_1''}) \otimes \dots \otimes 
\widehat {(\del^{(1)}_{a'_i}\otimes \del^{(2)}_{a''_i})} \otimes \dots
\otimes  (\del^{(1)}_{a_n'}\otimes \del^{(2)}_{a_n''})
\otimes (\del^{(1)}_{a'}\otimes \del^{(2)}_{a''_i}))\nn\\
&+ \sum_{a''\in A^{(2)}}\sum_{i=1}^n d^{(2)}_{a''_ia''} 
(Y^{(1)}\otimes Y^{(2)})(\Delta_{\mbn})\nn\\
&\qquad ((\del^{(1)}_{a_1'}\otimes \del^{(2)}_{a_1''}) \otimes \dots \otimes 
\widehat {(\del^{(1)}_{a'_i}\otimes \del^{(2)}_{a''_i})} \otimes \dots
\otimes  (\del^{(1)}_{a_n'}\otimes \del^{(2)}_{a_n''})
\otimes (\del^{(1)}_{a'_i}\otimes \del^{(2)}_{a''}))\nn\hskip-2cm\\
&- n\, d\, (Y^{(1)}\otimes Y^{(2)})(\Delta_{\mbn})
((\del^{(1)}_{a_1'}\otimes \del^{(2)}_{a_1''}) \otimes \dots
\otimes  (\del^{(1)}_{a_n'}\otimes \del^{(2)}_{a_n''}))\Big)\nn\\
\overset{(*)}{=}& \sum_{n\geq 3} \frac{1}{n!} x^{a_n'a_n''}\cdots
x^{a_1'a_1''}\Big((D^{(1)}+D^{(2)}-d) (Y^{(1)}\otimes Y^{(2)})(\Delta_{\mbn})
\nn\\
&\qquad ((\del^{(1)}_{a_1'}\otimes \del^{(2)}_{a_1''}) \otimes
\dots \otimes (\del^{(1)}_{a_n'}\otimes \del^{(2)}_{a_n''})) \hskip-2cm\nn\\ 
&- \sum_{a' \in A^{(1)}} r^{(1) a'} (Y^{(1)}\otimes Y^{(2)})
((\pi^*,id)(\Delta_{\mbn}))
((\del^{(1)}_{a_1'}\otimes \del^{(2)}_{a_1''}) \otimes 
\dots \otimes  (\del^{(1)}_{a_n'}\otimes \del^{(2)}_{a_n''})\otimes \del^{(1)}_{a'})
\hskip -5cm\nn\\
&- \sum_{a''\in A^{(2)}} r^{(2) a''} 
(Y^{(1)}\otimes Y^{(2)})((id,\pi^*)(\Delta_{\mbn}))
((\del^{(1)}_{a_1'}\otimes \del^{(2)}_{a_1''}) \otimes 
\dots \otimes  (\del^{(1)}_{a_n'}\otimes \del^{(2)}_{a_n''})\otimes
\del^{(2)}_{a''})\Big)\hskip -5cm\nn\\
=& \sum_{n\geq 3} \frac{1}{n!} \Big((D+d)(Y^{(1)}\otimes
Y^{(2)})(\Delta_{\mbn})(\gamma^{\otimes n}) \nn\\
&- \sum_{a' \in A^{(1)}} r^{(1) a'} (Y^{(1)}\otimes
Y^{(2)})((\pi^*,id)(\Delta_{\mbn}))(\gamma^{\otimes n}\otimes
\del^{(1)}_{a'})\nn\\
&- \sum_{a'' \in A^{(2)}} r^{(2) a''} (Y^{(1)}\otimes
Y^{(2)})((id,\pi^*)(\Delta_{\mbn}))(\gamma^{\otimes n} \otimes
\del^{(2)}_{a''})\Big). 
\label{t3.12}
\end{align}
To obtain $(*)$ write $\Delta_{\mbn}= \sum [\tau]
g^{[\tau] [\sigma]}\otimes [\sigma]$ as in \ref{ts3.1.1}
and apply Proposition \ref{ts2.6.3} to both tensor
factors of each summand. Furthermore, notice that the 
$[\tau], [\sigma]$ are homogeneous 
and $g^{[\sigma] [\tau]} =0$
unless 
$|E_{\sigma}|+|E_{\tau}|=n-3$ (\ref{t3.2}).

\med

On the other hand, applying Proposition \ref{ts2.7.2}, we obtain up to
quadratic terms
\begin{align}
E_0\Phi(\gamma) =&  \sum_{n\geq 3} \frac{1}{(n-1)!} 
\Big(\sum_{a' \in A^{(1)}} r^{(1) a'} (Y^{(1)}\otimes Y^{(2)})(\Delta_{\mbn})(\gamma^{\otimes
n-1}\otimes \del^{(1)}_{a'}\otimes
\del^{(2)}_0)\nn\\
&\quad + \sum_{a'' \in A^{(2)}} r^{(2) a''}
(Y^{(1)}\otimes Y^{(2)})(\Delta_{\mbn})(\gamma^{\otimes n-1} \otimes \del^{(1)}_0
\otimes
\del^{(2)}_{a''})\Big)\nn\\
=&
\sum_{n\geq 3} \frac{1}{n!} 
\Big(\sum_{a' \in A^{(1)}} r^{(1) a'} (Y^{(1)}\otimes Y^{(2)})(id,\pi_*)
(\Delta_{\overline{M}_{n+1}})
(\gamma^{\otimes n}\otimes \del^{(1)}_{a'})\nn\\
&\quad +  \sum_{a'' \in A^{(2)}} (Y^{(1)}\otimes Y^{(2)})
(\pi_*,id)(\Delta_{0\overline{M}_{n+1}})
(\gamma^{\otimes n} \otimes \del^{(2)}_{a''})\Big)
\label{t3.13}
\end{align}

Applying the formula (\ref{t3.8}), we see that the sum of 
(\ref{t3.12}) and (\ref{t3.13}) is just the the quasi--homogeneity
condition for $E$ and therefore $E$ is an Euler field.

\med

{\it Identity.}

\med

The proposed identity $\del^{(1)}_0 \otimes \del^{(2)}_0$ is a flat field
by definition.
Furthermore,
\begin{multline}
Y_3(\del^{(1)}_{a'}\otimes\del^{(2)}_{a''} \otimes
\del^{(1)}_{b'}\otimes\del^{(2)}_{b''}
\otimes \del^{(1)}_0 \otimes \del^{(2)}_0) =\\ 
Y^{(1)}_3(\del^{(1)}_{a'}\otimes\del^{(1)}_{b'} \otimes \del^{(1)}_0) 
Y^{(2)}_3(\del^{(2)}_{a''}\otimes\del^{(2)}_{b''} \otimes \del^{(2)}_0)
= g_{a'b',a''b''}
\label{t3.14}
\end{multline}
and for $n \geq3$ by Proposition \ref{ts2.7.2} and (\ref{t3.9})
\begin{align}
Y_n&((\del^{(1)}_{a_1'}\otimes\del^{(2)}_{a_1''})\otimes \dots \otimes
(\del^{(1)}_{a_{n-1}'} \otimes\del^{(2)}_{a_{n-1}''}) \otimes 
(\del^{(1)}_0 \otimes \del^{(2)}_0))\nn\\
& =(Y^{(1)}\otimes Y^{(2)})((\pi_*,\pi_*)\Delta_{\mbn})
((\del^{(1)}_{a_1'}\otimes\del^{(2)}_{a_1''})
\otimes \dots \otimes (\del^{(1)}_{a_{n-1}'}
\otimes\del^{(2)}_{a_{n-1}''}))\nn\\
&=0 \label{t3.15}
\end{align}
which proves that $\del^{(1)}_0 \otimes \del^{(2)}_0$ is indeed an
identity by Remark \ref{ts2.7.1}.
The weight of this identity can be read off the Euler field as
$d+d-d=d$,
proving the theorem.

\med

\subsubsection{\label{ts4.2.3}Remarks}

\sm
The condition that the weights of the identities are equal can be
met by a rescaling of the Euler fields as long as 
not only one of the weights is $0$. In the following, we will always
assume this when considering the tensor product.

\sm
Since, given a metric and the multiplication on the fibers of a
Frobenius manifold,
the identity is uniquely determined ---cf\ .[M1]---, the above
identity is the only identity
compatible with the choice of the tensor metric (\ref{t3.10}). 

\sm
The theorem, however, contains no such uniqueness property for the
Euler field,
but there are several reasons for the choice of this particular
type of Euler field. If the $E_1$--part is regarded as providing the
operator $\Cal V$ of (\ref{t1.8}), then our choice of $E_1$ for the tensor product
is equivalent up to the shift by $d$ which is necessary to accommodate
the dependence of the tensor product on the diagonal
in $H^*(\mbn \times \mbn)$ to the natural definition:
\begin{equation}
{\Cal V} := {\Cal V}^{(1)} \otimes {\mathrm id} + {\mathrm id} \otimes {\Cal V}^{(2)}.
\label{t3.16}
\end{equation}

As remarked in [M1], 
if the action of ${\mathrm ad}(E)$ is semi--simple on $H$, 
there is a natural grading of $H$ induced by
the action of
${{\mathrm ad}}(E)$, shifted by $d_0$. This grading basically
fixes
the $E_1$ component. In the setting of quantum cohomology, this
grading is just (half) the usual grading for the cohomology groups. The
additivity is just the
fact that under the K\"unneth formula 
the total degree of a class is the sum of the degrees of the two
components.
The natural grading on the space of  $H^{(1)}\otimes H^{(2)}$
is consequently 
given by the grading operator
${{\mathrm ad}}(E^{(1)} \otimes {\mathrm id} + {\mathrm id} \otimes E^{(2)})$ shifted by $d$,
so
that the
tensor product of $\del^{(1)}_a$ and $\del^{(2)}_b$ of
degrees 
$\delta^{(1)}_{a'}$ and $\delta^{(2)}_{a''}$ is of degree
$\delta^{(1)}_{a'} +d+\delta^{(2)}_{a''}+d-d$. Recalling that $d_a$
was the eigenvalue
of  ${{\mathrm -ad}}(E)$, we 
obtain $d_{a'a''}= d^{(1)}_{a'} + d^{(2)}_{a''} -d$.

\sm

In the physical realm of topological field theories [DVV],
the above argument
for the choice of $E_1$ just reflects the additivity of a $U(1)$
charge.

\med

The choice for $E_0$ is motivated by quantum
cohomology where
the $E_0$--part corresponds to
the canonical class. Thus, the definition of 
$E_0 = E^{(1)}_0 \otimes \del^{(2)}_0 +\del^{(1)}_0 \otimes E^{(2)}_0$ 
corresponds to the formula 
$K_{X \times Y} = K_X\otimes 1 + 1 \otimes K_Y$.
More generally, it corresponds to the map
$H^*(V) \times H^*(W) \ra H^*(V\times W):(v,w) \ra pr_1^*(v) + pr_2^*(w)$
which generally reflects the structure of the tensor product
of Frobenius manifolds in the presence of flat identities, see
Section \ref{idtenfrob} below.

\sm

Furthermore, in view of (\ref{t3.12}) and Lemma \ref{ts3.2.2},
$E_0$ seems to be the only possible choice, if one 
postulates (\ref{t3.16}).

\subsection{\label{ts4.2}The tensor product for two germs of pointed
Frobenius manifolds}

\bigskip

Due to the following main Theorem of this section, we can define
the tensor product in the category of pointed germs of Frobenius manifolds.

\med

\subsubsection{Theorem}
\label{convgerm}

{\it
The tensor product potential of two convergent potentials is convergent.}
\med

\subsubsection{\label{ts4.2.1}Definition}
Given two germs of pointed Frobenius manifolds $({\bold M}^{(1)},m^{(1)}_0)$ 
and $({\bold M}^{(2)}, m^{(2)}_0)$, let $(H^{(1)},g^{(1)},\Phi^{(1)})$ 
and $(H^{(2)},g^{(2)},\Phi^{(2)})$ be the associated formal Frobenius 
manifolds.
We define the tensor product $({\bold M},m_0)$ 
of $({\bold M}^{(1)}, m^{(1)}_0)$ 
and $({\bold M}^{(2)}, m^{(2)}_0)$ to be the associated germ of a pointed 
Frobenius manifold.

\bigskip

We will now prove the main Theorem of this section in several steps
starting with invertible 1-dimensional CohFTs and proceeding to
full generality.
In the course of the proof, we will utilize a Theorem on complex series 
cited below for convenience
(cf.\ e.g.\ Grauert, Einf\"uhrung in die Funktionentheorie mehrerer 
Ver\"anderlicher, Satz 1.1).

\med

\subsubsection{Theorem.\label{gen}}
{\it
Let ${\bold z}_1 \in  \stackrel {{\rm o}} {{\Bbb C}^n}:= \{{\bold z}
= (z_1, \dots , z_n) \in {\Bbb C}^n \bigm| z_k \neq 0; 1 \leq k \leq n\}$.
If the power series 
$\sum_{\boldsymbol{\nu} = 0}^{\infty} a_{\boldsymbol{\nu}} 
{\bold z}^{\boldsymbol{\nu}}$ converges at ${\bold z_1}$, 
then the series is uniformly convergent
inside the polycylinder $P_{\bold z_1} := \{{\bold z} \in {\Bbb C}^n 
\bigm| |z_k| < |z_k^{(1)}|\}$.} 

\med 
\subsubsection{Invertible rank one CohFTs}

\sm

For invertible rank one CohFTs (i.e.\ $C_3\neq 0$) we have the following 
property:
\subsubsection{\label{kmzs3.4.2}Theorem (3.4.2 of [KMZ])}  
{\it
Define the  bijections
\begin{equation} 
\NIC\LR\frac{x^3}6+x^4\,k[[x]]\LR 1+\eta\,k[[\eta]] \,,\label{kmz3.4}
\end{equation}
where the first map assigns to a theory $\A$ its 
potential $\Phi_\A(x)$ and the second map is defined by
\begin{equation}
\Phi(x)\,\leftrightarrow\, U(\eta)=\int_0^\infty e^{-\Phi''(\eta x)/\eta}\,dx
\label{kmz3.5}
\end{equation}
or alternatively by assigning 
to $\Phi(x)=\frac16x^3+\dots$ 
the power series $U(\eta)=\sum_{n=0}^\infty B_n\eta^n$
where $x=\sum B_n\frac{y^{n+1}}{(n+1)!}=y+\cdots$ 
is the inverse power series of $y=\Phi^{(2)}(x)=x+\cdots$. 
Then the tensor product of 1--dimensional CohFTs corresponds to 
multiplication 
in $1+\eta k[[\eta]]\,$:
$U_{\A^{(1)}\otimes\A^{(2)}}(\eta)=U_{\A^{(1)}}(\eta)\,U_{\A^{(2)}}(\eta)$. 
The coefficients of $-\log U_\A(\eta)$ are
the canonical coordinates of $\A\,$.
}

\sm

\subsubsection{\label{3.5}Explicit formulas.} 
The above theorem can be used to give explicit formulas
for the coefficients of $U(\eta)$ in terms of the coefficients of $\Phi(x)$;
see [KMZ] section 3.5.
The explicit law for the tensor product of two normalized invertible CohFTs 
in terms of the coefficients of their
potential functions can be derived by combining these formulas with the 
identity $U_{\A^{(1)}\otimes\A^{(2)}}(\eta)=U_{\A^{(1)}}(\eta)\,U_{\A^{(2)}}(\eta)$:
\begin{align*} 
C_4 &=C_4^{(1)}+C_4^{(2)}\,,\\  C_5 &=C_5^{(1)}+5C_4^{(1)}C_4^{(2)}+C_5^{(2)}\,,\\ 
 C_6 &=C_6^{(1)}+(8\,{C_4^{(1)}}^2+C_5^{(1)})\,C_4^{(2)}+C_4^{(1)}\,(8\,{C_4^{(2)}}^2+C_5^{(2)})+C_6^{(2)}\,,\\
 C_7 &=C_7^{(1)}+ (35\,C_4^{(1)}\,C_5^{(1)}+14\,C_6^{(1)})\,C_4^{(2)}+(61\,{C_4^{(1)}}^2\,{C_4^{(2)}}^2+33\,
{C_4^{(1)}}^2C_5^{(2)}+ \\ 
&\qquad 33\,C_5^{(1)}\,{C_4^{(2)}}^2 
 +19\,C_5^{(1)}\,C_5^{(2)}) + 
C_4^{(1)}\,(35\,C_4^{(2)}\,C_5^{(2)}+14\,C_6^{(2)}) +C_7^{(2)}\,,
\qquad\ldots 
\end{align*}
\med

\subsubsection{Proposition.\label{inv}} 
{\it
The tensor product potential of two convergent invertible rank one CohFTs
is a convergent invertible rank one CohFT.}

\sm

{\bf Proof.}
First assume that $C^{(1)}_3=C^{(2)}_3=1$.

\sm

As in Theorem \ref{kmzs3.4.2} write 
the inverse power series of the second derivatives of the potentials
$$\Psi^{(1)}:= \sum \frac{B^{(1)}_n}{(n+1)!} y^{n+1} \qquad 
\Psi^{(2)} := \sum \frac{B^{(2)}_n}{(n+1)!} y^{n+1}.$$
These series are convergent, if the respective potentials are and so is the 
product of their positive counterparts 
$|\Psi^{(1)}|:= \sum \frac{|B^{(1)}_n|}{(n+1)!} y^{n+1}$ and  
$|\Psi^{(2)} |:= \sum \frac{|B^{(2)}_n|}{(n+1)!} y^{n+1}$ as well as its derivative:

$$\widehat{\Psi}:= \frac{\del}{\del y} (|\Psi^{(1)}|\,|\Psi^{(2)}|) 
= \sum_{n \geq 0}  \sum_{i=0}^n \frac{n+1}{i+1}
\frac{|B^{(1)}_i B^{(2)}_{n-i}|}{i!(n+1-i)!} y^{n+1}.$$ 

By [KMZ] (see Theorem \ref{kmzs3.4.2} above) 
we have the following expansion for the inverse of the second 
derivative of the tensor potential $\Psi$:
\begin{equation}
\Psi = \sum_{n \geq 0} \sum_{i=0}^n
\frac{B^{(1)}_i B^{(2)}_{n-i}}{(n+1)!} y^{n+1}. 
\end{equation}
This series is dominated by $\widehat{\Psi}$ at any point inside the
domain of convergence of $\widehat{\Psi}$, since 
$$\frac{n+2}{i+1}\binom{n+1}{i} > 1, \text{ for } 0\leq i \leq n,$$
and thus
$$
|\sum_{i=0}^n \frac{B^{(1)}_i B^{(2)}_{n-i}}{(n+1)!}| \leq
\sum_{i=0}^n \frac{|B^{(1)}_i B^{(2)}_{n-i}|}{(n+1)!} <
\sum_{i=0}^n \frac{n+1}{i+1} \frac{|B^{(1)}_i B^{(2)}_{n-i}|}{i!(n+1-i)!}
$$
proving the convergence of the tensor product potential if $C^{(1)}_3=C^{(2)}_3=1$.

\sm

In case that $C^{(1)}_3$ and/or $C^{(2)}_3$ are not equal to one, but not equal to zero,
we can scale $\del^{(1)}$ and $\del^{(2)}$ in such a way that they are one.
Notice the following scaling behavior for potentials of
one--dimensional theories
$$
\Phi_{\lambda \del}(x) = \Phi_{\del}(\lambda x)=\Phi_{\del}^{\lambda}(x),
$$ 
where the subscript $\del$ refers to the chosen basis vector of the 
theory and $\Phi_{\del}^{\lambda} (x)$ is the potential for the scaled 
theory $C_i^{\lambda} = \lambda^i C_i, i\geq3$ which is  
convergent, if $\Phi$ is. 

Furthermore, for the tensor product potential
of two such theories: 
\begin{multline*}
\Phi^{(1)}_{\del^{(1)}}\otimes \Phi^{(2)}_{\del^{(2)}}(\lambda \mu x)=
\Phi_{\del^{(1)}\otimes \del^{(2)}} (\lambda \mu x)=
\Phi_{\lambda \del^{(1)}\otimes \mu \del^{(2)}} (x)\\
=\Phi^{(1)}_{\lambda \del^{(1)}}\otimes \Phi^{(2)}_{\mu \del^{(2)}}(x)=
\Phi^{(1) \lambda}_{\del^{(1)}}\otimes \Phi^{(2) \mu}_{\del^{(2)}}(x).
\end{multline*}

Choosing the appropriate scalings the right hand side is the tensor product 
of two convergent potentials 
with $C_3^{(1)}$=$C_3^{(2)}$=$1$ which converges by 
the first argument. It follows that $\Phi^{(1)}_{\del^{(1)}}\otimes \Phi^{(2)}_{\del^{(2)}}$
also converges, proving the proposition.

\med

\subsubsection{The case of general rank one CohFTs}

\sm

In the previous section, we used the fact that we have a good
handle on the tensor product potential in the case that the 
two rank one theories are invertible. We will show below
that one can basically use the same formula even if the 
theories in question are not--necessarily invertible.

\sm

Denote by $C_n$ the coefficients of the tensor product 
potential of two rank one CohFTs. 

\sm

For a monomial 
$p = \text{const.} \times C_{i_1} \dots C_{i_n}$ define the degree
$deg(p):= i_1 + \cdots + i_n$ and the length
$length(p):=n$

\med

\subsubsection{Lemma.}
{\it In the notation of the previous section
\begin{equation}
C_n = Y^{(1)}\otimes Y^{(2)} (\Delta_{0n})(\del^{\otimes n})
= P_n (C^{(1)}_3, \dots, C^{(1)}_n, C^{(2)}_3, \dots, C^{(2)}_n)
\end{equation}
where $P_n(C^{(1)}_3, \dots, C^{(1)}_n, C^{(2)}_3, \dots, C^{(2)}_n)$ is a universal polynomial.
Furthermore,} 
\begin{equation}
P_n =  \sum \text{monomials $p_{n,(k^{(1)},k^{(2)}),(l^{(1)},l^{(2)})}^{(i)}$ 
in the $C^{(1)}_i, C^{(2)}_j$}
\end{equation}
{\it 
with $bideg(p_{n,(k^{(1)},k^{(2)}),(l^{(1)},l^{(2)})}^{(i)})= (k^{(1)},k^{(2)})$ , $ bilength(p_{n,(k^{(1)},k^{(2)}),(l^{(1)},l^{(2)})}^{(i)}) = (l^{(1)},l^{(2)})$
and the bi--degrees and bi--lengths satisfy:
\begin{equation}
k^{(1)}-2l^{(1)} = k^{(2)} -2l^{(2)}= n - 2\;
\text { and }\; l^{(1)}+l^{(2)}= n-1.
\label{res}
\end{equation}
}

\sm

{\bf Proof.} 
Just express $Y^{(1)}\otimes Y^{(2)} (\Delta_{0n})(\del^{\otimes n})$
as a sum over trees $(\tau^{(1)},\tau^{(2)})$. 
The restriction then follows from the
observation that $bideg(p)$ is $(|F_{\tau^{(1)}}|,|F_{\tau^{(2)}}|)$ and 
$bilength(p)$ is $(|V_{\tau^{(1)}}|,|V_{\tau^{(2)}}|)$.
Finally, notice that $|E_{\tau^{(1)}}|+|E_{\tau^{(2)}}| = n-3$.

\med

\subsubsection{Lemma.\label{resl}}
{\it For a fixed $n$ there is a unique way of extending
a monomial $p_{red}^n$ in $C^{(1)}_i,C^{(2)}_j; 4 \leq i,j  \leq n$ of given
bi--degree $(k^{(1)},k^{(2)})$ and bi--length $(l^{(1)},l^{(2)})$ into a
monomial $p^n$ in the $C^{(1)}_i,C^{(2)}_j; 3 \leq i,j \leq n$, s.t.
the bi--degree and bi--length of $p^n$ satisfy
the equations (\ref{res}) and $p^n$ coincides
with $p^n_{red}$ for $C^{(1)}_3=C^{(2)}_3=1$.}

\sm

{\bf Proof.}
The monomials must be of the form 
$p^n = p_{red}^nC_3^{(1) i}C_{3}^{(2) j}$ and
using the restrictions (\ref{res}) we find:
\begin{equation}
i= n-2-k^{(1)}+2l^{(1)} \; \text { and } \; j= n-2-k^{(2)}+2l^{(2)}.
\label{pow}
\end{equation}

\med

\subsubsection {Corollary.\label{ext}}
{\it
The universal polynomials $P_n$ are given by the unique polynomials
extending the $C_n =: p^n_{red}$ given in [KMZ].}

\med

\subsubsection{Proposition.\label{noninv}}
{\it
The tensor product of two convergent rank one CohFTs is again convergent.
}

\sm

{\bf Proof.}
Given the potentials $\Phi^{(1)}$ and $\Phi^{(2)}$, we assume after scaling 
that $C^{(1)}_3, C^{(2)}_3 \in \{0,1\}$. Denote by
$\widehat {\Phi}^{(1) / (2)}$ the potential
with $\widehat{C_3}^{(1) / (2)}=1$ 
and $\widehat{C_i}^{(1) / (2)}= 
C_i^{(1) / (2)}, i\geq 4$.
These potentials are both convergent and invertible. Using the proposition 
for convergent and invertible potentials, we obtain that their 
tensor potential $\widehat{\Phi}$ is convergent. 
Now, due to the Corollary \ref{ext}   
there is a unique power series $\widetilde{\Phi}\in {\Bbb C}[[C_3^{(1)},C_3^{(2)},x]]$
extending $\widehat{\Phi}$, s.t. 
$\widetilde{\Phi}|_{C^{(1)}_3=C^{(2)}_3=1} = \widehat{\Phi}$ and the conditions
(\ref{res}) are satisfied.

\sm

First, assume that only one of the potentials is not invertible 
say  $C^{(1)}_3=1, C_3^{(2)}=0$. Regarding the power series 
$\widetilde{\Phi}|_{C^{(1)}_3=1}=:\widetilde{\Phi}_1 \in {\Bbb C}[[C_3^{(2)},x]]$,
notice that
$\widetilde{\Phi}_1$ converges at all points $(1,x_0)$ with
$x_0$ inside the domain of convergence of $\widehat{\Phi}$ and is therefore
---again by Theorem \ref{gen}---
convergent at points $(0,x_0)$. However, 
$\widetilde{\Phi_1}|_{C^{(2)}_3=0} = \Phi$ and thus $\Phi$ is also
convergent.

\sm

In case that both $C^{(1)}_3=0$ and $C_3^{(2)}=0$, we see that 
$\widetilde{\Phi}|_{C^{(1)}_3=1,C^{(2)}_3=1} = \widehat{\Phi}$ and 
$\widetilde{\Phi}$ converges at all points $(1,1,x_0)$ with
$x_0$ inside the domain of convergence of $\widehat{\Phi}$. Therefore
---again by Theorem \ref{gen}---
it is also convergent at points $(0,0,x_0)$. Now, 
$\widetilde{\Phi}|_{C^{(1)}_3=C^{(2)}_3=0} = \Phi$ and we again obtain that
$\Phi$ is convergent.

\med

\subsubsection{Proposition}
\label{zerotensor}
{\it
In case that both potentials are non-invertible, i.e.\ $C^{(1)}_3=C^{(2)}_3=0$,
we even have that $\widetilde{\Phi}|_{C^{(1)}_3=C^{(2)}_3=0}= \Phi \equiv 0$.}

\sm

{\bf Proof.}
By Lemma \ref{resl}, all summands of $\Phi$ are of the form
$p^n = p_{red}^nC_3^{(1)  i}C_{3}^{(2) j}$ with $p_{red}$ of given
bi--length $(l^{(1)},l^{(2)})$
and bi--degree $(k^{(1)},k^{(2)})$ and $i,j$ given by (\ref{pow}).
Furthermore from the last equation in (\ref{res}) we obtain 
$l^{(1)}+l^{(2)}=n-1-(i+j)$. 
Thus, using the inequalities
$$
k^{(1)} \geq 4 l^{(1)} \quad k^{(2)} \geq 4 l^{(2)}
$$
we find:
$$
0= 2n-4 -(k^{(1)}+k^{(2)}) +2(l^{(1)}+l^{(2)})- (i+j) \leq -2 + (i+j).
$$
So that $(i+j)\geq2$ and all $p^n$ vanish for $C^{(1)}_3=C^{(2)}_3=0$.

\med

\subsubsection{The higher dimensional case}

\sm

Given two formal Frobenius manifolds $(V^{(1)}$, $g^{(1)}$, $\Phi^{(1)})$ 
and $(V^{(2)}$, $g^{(2)}$, $\Phi^{(2)})$
with convergent potentials $\Phi^{(1)}$ and $\Phi^{(2)}$, denote the corresponding
ACFs by $Y^{(1)}$ and $Y^{(2)}$.

\sm

Let $(V,g,\Phi)$ be the tensor product formal Frobenius manifold.
Choosing a basis $(\del^{(1)}_a), a\in A$ resp.\ $(\del^{(2)}_b), b \in B$ 
for $V^{(1)}$ resp.\ $V^{(2)}$, the tensor potential in the tensor basis 
$\del_{ab} := \del^{(1)}_a \otimes \del^{(2)}_b$ takes the form
\begin{equation}
\Phi = \sum_{n = 3}^{\infty} \frac{1}{n!}
\sum_{(a_1,\dots,a_n)\atop(b_1,\dots,b_n)} x_{a_nb_n} \dots x_{a_1b_1}
Y^{(1)}\otimes Y^{(2)}(\Delta_{0n}) (\del_{a_1b_1}\otimes \dots \otimes \del_{a_nb_n})
\end{equation} 

\med

\subsubsection{Pure even case}

\sm

In the pure even case, we can consider
the points $y_{diag}$ whose coordinates are given by
$x_{ab} \equiv y$, $y \in {\Bbb C}$ constant 
$\forall a\in A$ and $b\in B$. The potential at these points reads
\begin{equation}
\Phi = \sum_{n = 3}^{\infty} \frac{1}{n!} y^n 
Y^{(1)}\otimes Y^{(2)}(\Delta_{0n})((\del^{(1)} \otimes \del^{(2)})^{\otimes n})
\end{equation}
with 
\begin{equation}
\del^{(1)} := \sum_{a\in A}\del^{(1)}_a \; \text{ and } \; 
\del^{(2)}:=\sum_{b \in B} \del^{(2)}_b.
\end{equation}

\med

\subsubsection{Proposition} 
{\it The potential of two convergent pure even
CohFTs is convergent.}

\sm

{\bf Proof.}
First scale the chosen basis in such a way that $|g^{ab}| \leq 1$.
Now, consider the series
$$
\sum_{n = 3}^{\infty} \frac{1}{n!} y^n 
|Y^{(1)}\otimes Y^{(2)}(\Delta_{0n})((\del^{(1)} \otimes \del^{(2)})^{\otimes n})|.
$$
Due to the condition $|g^{ab}| \leq 1$, this series is dominated by the tensor 
potential for two rank one CohFTs given by the coordinates:
\begin{equation}
C^{(1)}_n := |Y^{(1)}_n(\del^{(1) \otimes n})|, \qquad 
C_n^{(2)} := |Y^{(2)}_n(\del^{(2)\otimes n})|.
\label{Ccoord}
\end{equation}

Since the two given potentials $\Phi^{(1)}$ and $\Phi^{(2)}$ are 
convergent, so are their restrictions to the line $x^{(1)}_a \equiv y$ resp.\
$x^{(2)}_b \equiv y$ as power series in ${\Bbb C}[[y]]$. Thus they are
also absolutely
convergent and the positive counterparts of these 
restrictions are just the rank one CohFT given by the 
coordinates (\ref{Ccoord}). 
The potential of the 
tensor product of two convergent CohFTs of rank one is convergent 
by Proposition \ref{noninv}. Therefore, we
have convergence of the tensor potential $\Phi$ of $\Phi^{(1)}$ and 
$\Phi^{(2)}$ at some points $y_{diag}$, by the remarks above. 
Using the theorem on complex series \ref{gen}, we find that
the potential $\Phi$ is indeed convergent.

\med

\subsubsection{The general case}

\sm

Consider the underlying ${\Bbb Z}_{2}$ graded space of the theory
$V= V_0 \oplus V_1$ with
a basis $(\del_{a_i}\bigm| i \in I_0)$ of $V_0$ and
$(\del_{a_i}\bigm| i \in I_1)$ for some subsets $I_0,I_1$ of a
set $I=I_0 \amalg I_1$ with an order $<$.
Again chose the basis in such a way that $|g^{ab}|\leq 1$. 
Denote the dual coordinates of $V_0$ by
$(x_i\bigm| i \in I_0)$ and those of $V_1$ by $(y_i\bigm| i \in I_1)$.

\sm

The potential can now be written as
\begin{align}
\Phi &= \sum_n \sum_{a_n > \cdots > a_1 | a_i \in I_1}
y_{a_n} \dots y_{a_1} \Phi_{a_1, \dots, a_n}\\
\Phi_{a_1, \dots, a_n}&:= \sum_m \frac{1}{m!} 
\sum_{(b_m, \cdots b_1) \in I_0^{\times n}}
x_{b_m} \cdots x_{b_1} Y_{n+m}(\del_{a_1} \otimes \dots \otimes \del_{a_n}
\otimes \del_{b_1} \otimes \cdots \otimes \del_{b_m})
\end{align}

The potential $\Phi$ is by definition convergent, if all 
the $\Phi_{a_1, \dots, a_n}$ are.

\med

\subsubsection{Proposition}
{\it If $\Phi$ is the tensor product of two convergent series $\Phi$ and 
$\Phi^{(2)}$, then
all the $\Phi_{a_1b_1, \dots, a_nb_n}$ with $(a_i,b_i)\in I_1$ are convergent 
with the usual notation
for the variables and index sets for tensor products.}

\sm

{\bf Proof.}
Consider the following auxiliary pure even series in 
the even variables $x_{ij}\bigm|i \in I^{(1)}, j\in I^{(2)}$
\begin{equation}
\Psi = \sum_n \sum_{(a_1, \dots, a_n) \in I^{(1)\times n} \atop
(b_1, \dots, b_n) \in I^{(2) \times n} }
\frac{1}{n!} x_{a_nb_n} \dots x_{a_1b_1} 
|Y_n(\del_{a_1b_1} \otimes \dots \otimes \del_{a_nb_n})|
\end{equation}

To prove convergence we again look at points $x_i \equiv y$ with $y\neq 0$.
In the tensor case, i.e.\ $Y_n = Y^{(1)}\otimes Y^{(2)}(\Delta_{0n})$
this series is dominated by the tensor product potential of the following
two one-dimensional series: 
\begin{eqnarray}
\Psi^{(1)}= \sum_n \frac{1}{n!}y^nC^{(1)}_n \quad &
C^{(1)}_n = |Y_n(\sum_{i\in I^{(1)}}\del_{a'_i}^{(1)\otimes n})|\\
\Psi^{(2)}= \sum_n \frac{1}{n!}y^nC^{(2)}_n \quad &
C^{(2)}_n = |Y_n(\sum_{i\in I^{(2)}}\del_{a''_i}^{(2)\otimes n})|
\end{eqnarray}
where $I^{(1)}$ and $I^{(2)}$ are the index sets of the two original
theories $\Phi^{(1)}$ and $\Phi^{(2)}$. Notice that the ${\Bbb Z}_2$--grading and the
order of the $\del_{a_i}$ is irrelevant, since we take the absolute values
of the correlators. Because the series
$\Psi^{(1)}$ and $\Psi^{(2)}$ are convergent, by the assumption of convergence
of $\Phi^{(1)}$ and $\Phi^{(2)}$, their tensor product potential
converges by Proposition \ref{noninv}. Therefore, $\Psi$ converges as well.

\sm

Finally notice that
\begin{multline}
\frac{\del}{\del x_{a_1b_1}}\cdots \frac{\del}{\del x_{a_nb_n}}\Psi
\bigm|_{x_{a_ib_i} =0 : (a_i,b_i)\in I_1}=\\
 \sum_m \frac{1}{m!} \hskip -5pt
\sum_{((c_1,d_1), \cdots, (c_m,d_m)) \in I_0^{\times n}}\hskip -20pt
x_{c_md_m} \cdots x_{c_1d_1} 
|Y_{n+m}(\del_{a_1b_1} \otimes \dots \otimes \del_{a_nb_n}
\otimes \del_{c_1d_1} \otimes \cdots \otimes \del_{c_md_m})|\label{ders}
\end{multline}
and these functions again converge for some points $y_{diag}$ with
$x_{ij}(y_{diag}) \equiv y\bigm|(i,j)\in I_0$.
This shows that $\Phi_{a_1b_1, \dots, a_nb_n}$ is absolutely convergent at some
points on $y_{diag}$ of the above type and thus
the proposition follows by Theorem \ref{gen}.

\sm

Collecting all results, we arrive at the general Theorem \ref{convgerm}.

\section{The tensor product for Frobenius manifolds}

\sm

In this section we will only deal with analytic Frobenius manifolds.
\sm

\subsection{The exterior product of two Frobenius manifolds}

\sm

Given two Frobenius manifolds ${\bf M}^{(1)},{\bf M}^{(2)}$ we can
consider the vector bundle $T_{M^{(1)}} \boxtimes T_{M^{(2)}}$ 
on $M^{(1)}\times M^{(2)}$
which we call the exterior product bundle.
Since we have an affine flat structure on both $M$ and $N$, we also have
such a structure on $M^{(1)}\times M^{(2)}$ given 
by $T_{M^{(1)}}^f \boxplus T_{M^{(2)}}^f$ 
on $T_{M^{(1)}\times M^{(2)}}$. In addition 
we have an affine flat structure on $T_M\boxtimes T_N$ in the sense that
\begin{equation}
T_{M^{(1)}}^f \otimes T_{M^{(2)}}^f \otimes {\cal O}_{M^{(1)}\times M^{(2)}} 
\cong T_{M^{(1)}} \boxtimes T_{M^{(2)}}
\end{equation}

\pagebreak

\subsection{The tensor product relative to a pair of base--points}

\med

\subsubsection{Reminder}
As explained in section \ref{ts4.1.4} given any pointed Frobenius manifold 
$({\bf M},p)$ there is an associated
convergent formal Frobenius manifold $(T_p,g,\Phi_p)$ given by the expansion
of the potential $\Phi$ of $M$ at $p$ in terms of the coordinates $(x_i)$ of 
a chosen basis of flat vector fields $(X_i)$. 

\med

\subsubsection{The tensor product of Frobenius manifolds 
relative to a pair of base--points}
Using the notion tensor product for formal
Frobenius manifolds in the context of pointed Frobenius manifolds, we arrive
at the following construction:
 
\sm

Given a pair of points $(p,q)$ in the product ${M^{(1)}}\times M^{(2)}$
there is an associated convergent series 
$\Phi_{pq} := \Phi^{M^{(1)}}_p \otimes \Phi^{M^{(2)}}_q$
in the dual coordinates of  
$T_{p,{M^{(1)}}} \otimes T_{q,M^{(2)}}$,
which defines a Frobenius manifold structure
on the domain of convergence of $\Phi_{pq}$ (cf.\ section \ref{ts4.1.5}).
We will denote the resulting pointed Frobenius manifold
by:
\begin{equation}
(T_{p,{M^{(1)}}} \otimes T_{q,M^{(2)}},0).
\end{equation}
Its germ corresponds to the convergent formal Frobenius
structure
\begin{equation}
(\bigoplus {\Bbb C} \del_a \otimes \del_b,g_{M^{(1)}}\otimes g_{M^{(2)}}, 
\Phi_p^{M^{(1)}} \otimes \Phi_q^{M^{(2)}}).
\end{equation}
\med

\subsection{Frobenius structures on the affine exterior bundle}

In the construction of the previous section we have defined 
over each point $(p,q)\in M^{(1)}\times M^{(2)}$ a pointed
Frobenius manifold on a 
neighborhood $V_{pq}$ the zero section $s$ of the fiber of 
the exterior product bundle.
Let $V$ be the union of all the $V_{pq}$.
$$
\xymatrix{
V\ar@{^{(}->}[r]&T_{M^{(1)}}\boxtimes T_{M^{(2)}}\ar[d]\\
&M^{(1)}\times M^{(2)}\ar[ul]^s
}
$$

\med
\subsubsection{Definition}
A bundle will be called {\em a bundle of pointed Frobenius manifolds}
if there exists a neighborhood $V$ of the zero section s.t.\ the 
intersection
of each fiber with this neighborhood is a Frobenius manifold.
Let $V_p$ is the intersection of $V$ with the fiber at $p$ 
then it is naturally a pointed
Frobenius manifold with base--point zero. 
\med
\subsubsection{Examples}
\label{ex}
\begin{itemize}
\item[1)]
By the previous remarks, the exterior
product bundle is actually a bundle of pointed Frobenius manifolds.
\item[2)]
Every tangent bundle of a Frobenius manifold is naturally
a bundle of pointed Frobenius manifolds. 
This can be shown in two equivalent ways.
Either one uses the pointed Frobenius manifold $(\Phi,p)$ to define
a potential near zero on the fiber over $p$ or one uses the 
affine connection associated to the flat structure connection of the
Frobenius manifold to define on each fiber $T_{p,N}$ the potential 
$\Psi_p(\xi):= \Phi(p')$ where $p$ is the
development of a path joining $p'$ and $\xi\in T_{p,N}$ is the point of 
the development into the fiber at $p$ of the point $p'$. The latter
construction is defined locally since the local holonomy groups vanish,
due to the flatness of the connection. 
\end{itemize}
\med

\subsubsection{Definition}
A flat affine connection on an exterior product bundle 
over the Cartesian product of two Frobenius manifolds ${\bf M}^{(1)}$ and
${\bf M}^{(2)}$ which is an
extension of the linear connection defined by the canonical flat structure
is called
{\em tensor product connection} if it respects the flat 
and the Frobenius structures, i.e.\ it satisfies the following
two conditions:
\begin{itemize}
\item[i)]
Let $\theta_{\tau}:T_{M^{(1)}}\boxplus T_{M^{(2)}} \to T_{M^{(1)}}\boxtimes T_{M^{(2)}}$
be the map 
corresponding to the tensorial $(GL_n \otimes GL_m,{\Bbb C}^{nm})$ 
1--form defined by the affine connection (cf.\ e.g.\ [KN]) 
then $\theta_{\tau}$ induces
a linear map 
\begin{equation}
\theta_{\tau}^f:T_{M^{(1)}}^f\boxplus T_{M^{(2)}}^f \to 
T_{M^{(1)}}^f\boxtimes T_{M^{(2)}}^f
\end{equation}
\item[ii)]
The parallel displacement w.r.t.\ the affine tensor projection
preserves the germs of Frobenius manifolds. I.e.\
for all local horizontal lifts $\tilde {x_t}$ a of curve $x_t$ in
$M^{(1)}\times M^{(2)}$ with $x_0 =(p,q)$ 
and $\tilde{x}_0= 0  \in T_{M^{(1)}}\boxtimes T_{M^{(2)}}|_{(p,q)}$
into $V \subset T_{M^{(1)}}\boxtimes T_{M^{(2)}}$, where $0$ is
the zero section of the corresponding vector bundle:
\begin{equation}
(T_{M^{(1)}}\boxtimes T_{M^{(2)}}|_{(p,q)},\tilde x_0) 
\stackrel{\theta_{\tau}}{\cong} 
(T_{M^{(1)}}\boxtimes T_{M^{(2)}}|_{x_t},\tilde x_t)
\quad \forall \tilde{x}_t\in V. 
\label{tencheck}
\end{equation}
\end{itemize}
\sm

\subsubsection{Remark}
Since the linear connection is flat and torsion free, 
such a tensor product connection locally identifies the germs of 
Frobenius manifolds in different fibers of the bundle of pointed 
Frobenius manifolds
via affine parallel displacement along arbitrary path connecting the base--points.

\subsubsection{Proposition}
{\it
\label{thetatau}
If the Frobenius manifolds ${\bf M}^{(1)}$ and ${\bf M}^{(2)}$ 
both carry flat identities
then the affine connection defined by the 1--forms:
\begin{equation}
\begin{gathered}
\theta_{\tau}^U: T_{M^{(1)}}\boxplus T_{M^{(2)}} |_U
\to T_{M^{(1)}}\boxtimes T_{M^{(2)}}|_U\\
\theta_{\tau}(\del^{(1)}_a)= \delta_{a0} \del_{a0}, \quad
\theta_{\tau}(\del^{(2)}_b)= \delta_{0b} \del_{0b}\nn
\end{gathered}
\end{equation}
is a tensor product connection.
Here again $(\del^{(1)}_a, \del^{(2)}_b)$  
is the restriction of a chosen basis 
of $T^f_{M^{(1)}}\boxplus T^f_{M^{(2)}} $ 
and $(\del_{ab} = \del^{(1)}_a,\otimes\del^{(2)}_b)$ the tensor basis of
$T^f_{M^{(1)}}\boxtimes T^f_{M^{(2)}}$
}
\sm

{\bf Proof.}
It is clear that the locally defined forms glue together and
that the condition i) is met.
The proof of the condition ii) is given below by calculating  
the respective ACFs. We give the proof including odd coordinates.

\med

\subsubsection{\label{ts4.3.1}Lemma}
{\it
Let $\{Y^p_n\}$ be the ACFs
corresponding
to $\Phi_p$  and
$\{Y^{p'}_n\}$  be the ACFs 
corresponding to $\Phi_p'$ where $p'$ is some point which lies
inside the domain of convergence of the potential $\Phi_p$.
Let $x^a(p')=:x_0^a$ be the
$x$--coordinates of this point.
The new operadic ACFs are then given by, see (\ref{t4.3}):

\sm
For any stable $n$--tree $\tau$:
\begin{multline}
Y^{p'}(\tau) (\del_{a_1} \otimes \dots \otimes \del_{a_n})=
\sum_{N \geq 0} \frac{1}{N!}\sum_{(b_1,\dots,b_N): b_i \in A}
\eps(b|a)\, x^{b_N}_0 \cdots x^{b_1}_0\\
Y^p(\pi^*_{\{n+1,\dots ,n+N\}}(\tau))
(\del_{a_1} \otimes \dots \otimes \del_{a_n} 
\otimes \del_{b_1} \otimes \dots \otimes \del_{b_N}).
\label{t4.5}
\end{multline}
}

\sm

{\bf Proof.}
Inserting (\ref{t4.3}) into the definition of $Y(\tau)$ (\ref{t2.11}), we see
that the correlation
functions having a pre--factor $x^{b_N}_0 \cdots x^{b_1}_0$ are
those belonging to trees with $N-n$ tails added in an arbitrary
fashion to
$\tau$. The sum over all of these trees is just 
$\pi^*_{\{n+1,\dots ,n+N\}}(\tau)$, whence the Lemma follows.

\med

In order to prove the Proposition, we also 
need a Lemma about the diagonal $\Delta_{\mbs}$ which extends 
Lemma \ref{ts3.2.2}.

\med

\subsubsection{\label{ts4.3.2}Lemma}
{\it
For any two disjoint subsets $S,T \subset \n$
\begin{equation}
(\pi^*_S,\pi^*_T)(\Delta_{\Mb_{0\n\setminus (S\cup T)}}) 
= (\pi_{T*},\pi_{S*})(\Delta_{\mbn}).
\label{t4.6}
\end{equation}
}

\sm

{\bf Proof.}
Writing $(\pi^*_S,\pi^*_T)$ as $(\pi^*_S, id) \circ (id,\pi^*_T)$,
we obtain, after repeated application of Lemma \ref{ts3.2.2} in an
appropriate version, that 
$$
(\pi^*_S,\pi^*_T) (\Delta_{\Mb_{0\n\setminus (S\cup T)}})= 
(\pi^*_S \circ \pi_{T*}, id)(\Delta_{\Mb_{0\n\setminus S}}).
$$
Since $\pi_S$ and $\pi_T$ commute if $T \cap S =
\emptyset$,
we can prove the equality (\ref{t4.6}) again by Lemma \ref{ts3.2.2}.

\med

{\bf Proof of the Proposition \ref{thetatau}.} 
\sm
 
For any point $(p,q)\in M^{(1)} \times M^{(2)}$ 
denote the domain of convergence of $\Phi^{(1)}_p$ by $U_p$ and the
domain of convergence of $\Phi^{(2)}_q$ by $U_q$ and
set $U_{(p,q)}= U_p \times U_q$. 
\sm

Choose a point $(p',q') \in U_{(p,q)}$.
Let $p'$ have the
coordinates $x^{(1)}_{a'}(p')=x^{(1) a'}_0$
and $q'$  the 
coordinates $x^{(2)}_{b''}(q')=x^{(2) b''}_0$.

\sm
 
Denote the ACFs corresponding to $\Phi^{M^{(1)}}_p$ by
$\{Y^p_n\}$ 
and denote the ACFs corresponding to $\Phi^{M^{(1)}}_{p'}$ by $\{Y^{p'}_n\}$. 
Likewise denote the ACFs corresponding to $\Phi^{M^{(2)}}_q$ 
by $\{Y^q_n\}$   
 and denote the
ACFs for the expansion $\Phi^{M^{(2)}}_{q'}$ by $\{Y^{q'}_n\}$.
\sm

The correlation functions of the formal 
tensor product potential $\Phi_{pq}$ are be given by
$\{Y^{pq}_n = Y^p \otimes Y^q (\Delta_{\mbn})\}$. 
The ACFs of the potential 
$\Phi_{p'q'}$ will be denoted $\{Y^{p'q'}_n\}$.
\sm
Finally let 
$\theta_{\tau}((p',q'))$ be the image of the affine development 
w.r.t.\ $\theta_{\tau}$ of $0\in T_{(p',q')}$
into $T_{M^{(1)},p} \otimes T_{M^{(2)},q}$ and suppose 
that $\theta_{\tau}((p',q'))\in V_{pq}$.
\sm

Denote
the ACFs corresponding to the expansion 
of the holomorphic function given by $\Phi_{pq}$
at $\theta_{\tau}((p',q'))$ by $\{Y^{\theta_{\tau}((p',q'))}_n\}$.

\sm
Since $\theta_{\tau}$ is linear, the
point $\theta_{\tau}((p',q')) \in T_{M^{(1)},p} \otimes T_{M^{(2)},q}$ 
has the coordinates
$$
x_{a'b''}(\theta_{\tau}((p',q')) )= \delta_{a'0} x^{(1) a'}_0 + \delta_{0b''} x^{(2) b''}_0.
$$

The equation (\ref{tencheck}) in terms of these correlation functions reads:
\begin{equation}
Y^{\theta_{\tau}((p',q'))}_n = Y^{p'q'}_n \qquad \forall n\geq3.
\label{acfcheck}
\end{equation}

Now:
\begin{multline}
Y^{\theta_{\tau}((p',q'))}_n (\del_{a'_1 a''_1} \otimes \dots \otimes
\del_{a'_n a''_n})\\
= \sum_{N \geq 0} \frac{1}{N!} \sum_{l=0}^N
\hskip -10pt \sum \Sb (b'_1,\dots,b'_l)|b' \in A'\\
(b''_1,\dots,b''_{N-l})|b'' \in A'' \endSb \hskip -10pt \binom {N}{l} 
\,\eps(b'0|a'a'')\eps(0b''|a'a'')\,
x^{(2) b''_{N-l}}_0 \cdots x^{(2) b''_{1}}_0
x^{(1) b'_l}_0 \cdots x^{(1) b'_1}_0\\
(Y^p \otimes Y^q)
(\Delta_{\Mb_{0,N+n}})
(\del_{a'_1 a''_1} \otimes \dots \otimes \del_{a'_n a''_n}
\otimes \del_{b'_10} \otimes \dots \del_{b'_l0} \otimes
\del_{0b''_1} 
\otimes \dots \otimes \del_{0b''_{N-l}})\nn
\end{multline}
\begin{multline}
\shoveleft{=  \sum_{N \geq 0} \frac{1}{N!}\sum_{l=0}^N 
\sum \Sb (b'_1,\dots,b'_l)|b' \in A'\\
(b''_1,\dots,b''_{N-l})|b'' \in A'' \endSb \binom {N}{l}
\,\eps(b'|a')\eps(b''|a'')\eps(b'|a'')\eps(b''|a')\,}\\
x^{(2) b''_{N-l}}_0 \cdots x^{(2) b''_{1}}_0
x^{(1) b'_l}_0 \cdots x^{(1) b'_1}_0 	
(Y^p \otimes Y^q)
((\pi_{\{n+l+1,\dots,n+N\} *},\pi_{\{n+1,\dots, n+l\} *})
(\Delta_{\Mb_{0,N+n}}))\\
(\del_{a'_1 a''_1} \otimes \dots \otimes \del_{a'_n a''_n}
\otimes \del_{b'_1} \otimes \dots \del_{b'_l} \otimes
\del_{b''_1} 
\otimes \dots \otimes \del_{b''_{N-l}})
\label{t4.7}
\end{multline}
due to Proposition \ref{ts2.7.2}. 

\sm

On the other hand, tensoring the ACFs
$\{Y^{p'}_n\}$ and $\{Y^{q'}_n\}$  
and then utilizing Lemma \ref{ts4.3.1}
yields:
\begin{multline}
Y^{p'q'}_n (\del_{a'_1 a''_1} \otimes \dots \otimes
\del_{a'_n a''_n})
:= (Y^{p'} \otimes Y^{q'})(\Delta_{\mbn})
(\del_{a'_1 a''_1} \otimes \dots \otimes \del_{a'_n a''_n})\\
= \sum_{N \geq 0} \sum_{l=0}^N 
\sum \Sb (b'_1,\dots,b'_l)|b' \in A'\\
(b''_1,\dots,b''_{N-l})|b'' \in A'' \endSb \frac {1}{N!(N-l)!}
\,\eps(b'|a')\eps(b''|a'')\eps(b'|a'')\eps(b''|a')\,\\
\times x^{(2) b''_{N-l}}_0 \cdots x^{(2) b''_{1}}_0
x^{(1) b'_l}_0 \cdots x^{(1) b'_1}_0 
(Y^p \otimes Y^q)
((\pi^*_{\{n+1,\dots ,n+l\}},\pi^*_{\{n+l+1,\dots ,n+N\}})
(\Delta_{\Mb_{0,N+n}}))\\
(\del_{a'_1 a''_1} \otimes \dots \otimes \del_{a'_n a''_n}
\otimes \del_{b'_1} \otimes \dots \del_{b'_l} \otimes
\del_{b''_1} 
\otimes \dots \otimes \del_{b''_{N-l}}).
\label{t4.8}
\end{multline}

\sm

Applying Lemma \ref{ts4.3.2} with
$S= \{n=1,\dots ,n+l\}, T=\{n+l+1,\dots,n+N\}$,
we see that (\ref{t4.7}) and (\ref{t4.8}) and thus the multiplications,
respectively the potentials modulo quadratic terms, coincide. 

\med

\subsubsection{Remark}
One should view the existence of a tensor product connection as an expression
of the independence of the choice of base--points for the operation of forming
the tensor product of pointed germs of Frobenius manifolds.
More precisely, consider two germs of pointed Frobenius manifolds 
realized as two small neighborhoods of zero on ${\Bbb C}^{n_i}$
where $n_i$ are the appropriate dimensions. In this case, the pointed
tensor product is just the germ on the fiber of the exterior product
bundle over zero which we can again realize as some small open neighborhood of
zero. This neighborhood then contains all nearby germs (e.g.\ the germs
of tensor product with base--points near zero) via the affine connection.
Hence all continuations of the initial germ over zero will likewise
be continuations of these germs.

\sm

\subsection{The tensor product for Frobenius manifolds}
\label{idtenfrob}

In this section we will give a way to patch together all germs 
on the exterior product bundle. More precisely, we will
construct a Frobenius manifold which contains a submanifold parameterizing
all these germs. In order to give a general construction we will
need the following technical assumption.

\sm
\subsubsection{Definition}
 A flat identity on Frobenius manifold ${\bf M}$ is called {\em factorizable}
if $M = \bar M \times {\Bbb C}$, where the factor 
${\Bbb C}$ is coordinatized by the identity.

\sm
\subsubsection{Definition\label{tendef}}
Consider a commutative diagram of the type
\begin{equation}
\xymatrix{
T_{M^{(1)}}\boxtimes T_{M^{(2)}} \ar[d]\ar[r]^-{\Theta}
&\tau^*(T_N)\ar[dl]\ar[r]^{\hat \tau}&T_N\ar[d]\\
M^{(1)}\times M^{(2)}\ar[rr]_{\tau}\ar[d]^{p}&&N\\
\bar M^{(1)} \times \bar M^{(2)}\times {\Bbb C} \ar@{^{(}->}[rru]^-{i}&&
}
\label{diag}
\end{equation}
\med
where $(M^{(1)},T_{M^{(1)}}^f,g_{M^{(1)}},\Phi^{M^{(1)}})$, 
$(M^{(2)}, T_{M^{(2)}}^f, g_{M^{(2)}}, \Phi^{M^{(2)}})$ and
$(N, T^f_N, g_N, \Phi^N)$
are 
Frobenius manifolds 
with factorizable flat identities, 
$\tau$ is an affine map, which factors through 
$p: M^{(1)}\times M^{(2)}= \bar M^{(1)} \times {\Bbb C}\times \bar M^{(2)}
\times {\Bbb C} \ra  \bar M^{(1)} \times \bar M^{(2)}\times {\Bbb C}$ 
which is given in some fixed choice of coordinates on the factors
${\Bbb C}$ by
$p(m_1,x,m_2,y)= (m_1,m_2,x+y)$, $i$ is an embedding of affine flat manifolds  
and $\Theta$ is an isomorphism of metric bundles
with affine flat structure 
between the pulled back
tangent bundle of $N$  
and the exterior
product bundle over $M^{(1)} \times M^{(2)}$. Where the statement that
$\Theta$ is an isomorphism of metric bundles
with affine flat structure means that it is an isomorphism of metric bundles
and $\tau^*T^f_N= \Theta(T^f_{M^{(1)}}\otimes T^f_{M^{(2)}})$.

\med
We will call such a diagram a {\em  tensor product diagram} and
${\bf N}$ a tensor product manifold for ${\bf M}^{(1)}$ and ${\bf M}^{(2)}$
if it additionally preserves the structure of bundles of Frobenius manifolds
i.e.\ it satisfies the following conditions:
\begin{itemize}
\item[i)] 
for all points $(p,q)\in M^{(1)}\times M^{(2)}$:
\begin{multline}
(\bigoplus {\Bbb C} \; \del_a \otimes \del_b,g_{M^{(1)}}\otimes g_{M^{(2)}},
\Phi_p^{M^{(1)}} \otimes \Phi_q^{M^{(2)}}) 
\stackrel{\hat\tau \circ \Theta}{\cong}\\
(\bigoplus {\Bbb C} \; (\hat\tau\circ\Theta (\del_a \otimes \del_b)),g_N, 
\Phi^N_{\tau((p,q))}),
\label{germcond}
\end{multline}
where $\hat\tau\circ\Theta (\del_a \otimes \del_b) \in T_{\tau((p,q)),N}$

\sm
In other words, at all points in the 
image of $\tau$, $\hat\tau\circ\Theta$ gives an isomorphism of pointed
germs of Frobenius manifolds defined in Example \ref{ex}:
\begin{equation}
(T_{p,M^{(1)}}\otimes T_{q,M^{(2)}},0) \stackrel{\hat\tau \circ\Theta}{\cong} 
(T_{\tau((p,q)),N},0)=
({\bf N}, \tau((p,q)))
\end{equation}
\item[ii)]
The affine connection defined on $T_{M^{(1)}}\boxtimes T_{M^{(2)}}$
by the pullback of the canonical affine connection on $T_N$ ---
defined by the flat structure on $T_N$ and the canonical 1--form 
$\theta_{can}$--- is the tensor product connection
$\theta_{\tau}$ of the Proposition \ref{thetatau}:
\begin{equation}
\tau^*(\theta_{can})= \Theta \circ\theta_{\tau}.
\end{equation}
as maps.
\end{itemize}

\med

\subsubsection{Remark}
Due to the condition i) for a tensor product diagram the condition
ii) for a tensor product connection is already satisfied by the
pulled back affine connection cf.\ Example \ref{ex}. The condition i) 
for a tensor product connection then 
forces that $\tau$ is affine in the affine coordinates of the
source and target spaces. 

\subsubsection{Definition}
\label{isodef}
Two tensor product diagrams 

\begin{equation}
\xymatrix{
T_{M^{(1)}}\boxtimes T_{M^{(2)}} \ar[d]\ar[r]^-{\Theta}
&\tau^*(T_N)\ar[dl]\ar[r]^{\hat \tau}&T_N\ar[d]\\
M^{(1)}\times M^{(2)}\ar[rr]_{\tau}\ar[d]^-{p}&&N\\
\bar M^{(1)} \times \bar M^{(2)}\times {\Bbb C} \ar@{^{(}->}[urr]^-{i}&&
}
\label{unimod}
\end{equation}
 
\begin{equation}
\xymatrix{
T_{M^{(1)}}\boxtimes T_{M^{(2)}} \ar[d]\ar[r]^-{\Theta'}
&\tau^{\prime *}(T_N)\ar[dl]\ar[r]^{\hat \tau'}&T_{N'}\ar[d]\\
M^{(1)}\times M^{(2)}\ar[rr]_{\tau'}\ar[d]^-{p'}&&N'\\
\bar M^{(1)} \times \bar M^{(2)}\times {\Bbb C} \ar@{^{(}->}[urr]^-{i'}&&\\
}
\label{otherdiag}
\end{equation}
are called {\em  equivalent}
if there exist open neighborhoods $U_N, U_{N'}$ of the images of 
$\tau, \tau'$
and an isomorphism $\phi$ of
of Frobenius manifolds $U_N\to U_{N'}$ s.t.\ the induced  diagram
satisfying the conditions
of a tensor product diagram is commutative.
\begin{equation}
\xymatrix{
T_N|_{U_N}\ar@(ur,ul)@{.>}[rrrr]^{\phi_*}\ar[dd]&\tau^*T_N\ar[l]_{\hat \tau}
\ar[dr]
&T_{M^{(1)}}\boxtimes T_{M^{(2)}}\ar[d]\ar[l]_{\Theta}\ar[r]^{\Theta'}&
\tau^{\prime *}T_{N'}\ar[dl]\ar[r]^{\hat\tau'}&T_{N'}|_{U_{N'}}\ar[dd]\\
&&M^{(1)}\times M^{(2)}\ar[lld]^{\tau}\ar[rrd]^{\tau'}
\ar@<-.3cm>[d]_p\ar@<.3cm>[d]^{p'}&&\\
U_N\ar@(dr,dl)@{.>}[rrrr]_{\phi}&&
\bar M^{(1)} \times \bar M^{(2)}\times {\Bbb C} \ar@{^{(}->}[rr]^-{i'}
\ar@{_{(}->}[ll]_-{i}&&U_{N'}\\
}
\label{iso}
\end{equation}

Note that we take the notion of
isomorphism of Frobenius manifold in the strict sense
that all data should be compatible and we do not allow for instance a conformal
change in the metric as in [D2].

\sm

\subsubsection{Theorem}
\label{exist}
{\it For any two  
Frobenius manifolds with factorizable flat identities there exists
a tensor product diagram and hence a tensor product manifold
of these two manifolds. Furthermore, any two tensor product diagrams for two
given Frobenius manifolds ${\bf M}^{(1)}$ and ${\bf M}^{(2)}$ are equivalent.}

\sm

{\bf Proof.}

\med
{\bf Construction.}
\sm
We construct a tensor product manifold and the structure isomorphism
for Frobenius manifolds
with factorizable flat identities (for a generalization see below). 
We start from a pointed cover which is a refinement of the cover given by
the domains of convergence of the various $\Phi^{(1)}_p$ and $\Phi^{(2)}_q$.
Another choice of refinement would lead to an 
equivalent diagram.
\med
We define 
${\Cal W} = \{U_{(p,q)}|(p,q) \in M^{(1)} \times M^{(2)}\}$
where the $U_{(p,q)}$ are defined as above. 
\sm
The affine structure provides the affine transition functions for
the respective coordinate maps $\varphi_U$:
$$
\varphi_U^V = \varphi_V \varphi_U^{-1}: 
\varphi_U(U\cap V) \to \varphi_V(V\cap U) 
\quad \forall \; U, V \in {\Cal W}, \; U \cap V \neq \emptyset
$$
which can be written in matrix form relative to the chosen basis
$$
\varphi_U^V =
\left( \begin{array}{cc}
A_U^V& \xi_U^V\\
0&1
\end{array}
\right)
$$
They satisfy the conditions
\begin{equation}
\left( \begin{array}{cc}
A_U^V& \xi_U^V\\
0&1
\end{array}
\right)=
\left( \begin{array}{cc}
A_V^U& -A_V^U\xi_V^U\\
0&1
\end{array}
\right),\quad
\left( \begin{array}{cc}
A_U^W& \xi_U^W\\
0&1
\end{array}
\right)
=
\left( \begin{array}{cc}
A_U^V A_V^W& A_V^W\xi_U^V+ \xi_V^W\\
0&1
\end{array}
\right)
\label{equi}
\end{equation}
where $\varphi_U^V$ are  product transition functions, i.e.\
$A_U^V = (A^{(1)}\oplus A^{(2)})_U^V$ and
$\xi_U^V = x_U^V + y_U^V$.

Notice that the $U_{p,q}$
can be decomposed as
$U_{(p,q)}= \check U_{(p,q)} \times {\Bbb C}$ 
where we choose the direction ${\Bbb C}$ 
to be the ``anti--diagonal'' in ${\Bbb C}^2$ given by $e^{(1)}-e^{(2)}$. 
Furthermore, the map $p$ has
a section $s$ defined by $s(m_1,m_2,x):=(m_1,\frac{1}{2}x, m_2,\frac{1}{2}x)$.

\sm
Using the tensor product connection we can map each $U$ into the fiber
of the exterior product bundle over its base--point, say $(p,q)$. This yields
an embedding of $\check U$ into $V_{pq}$, where $V_{p_iq_i}$ is, as 
defined above, the domain of convergence
of $\Phi^{(1)}_{p_i} \otimes \Phi^{(2)}_{q_i}$.
Let $\theta_{\tau}$ be the tensor connection of \ref{thetatau}.
Denote the affine parallel displacement w.r.t.\ $\theta_{\tau}$
of $T_{M^{(1)}}\boxtimes T_{M^{(2)}}|_{(p,q)}$
into $T_{M^{(1)}}\boxtimes T_{M^{(2)}}|_{(p_0,q_0)}$
by $(\theta_{\tau})^{(p_0,q_0)}_{(p,q)}$. We define a
pointed coordinate neighborhood,
to be a pair $(U, (p_0,q_0))$ 
s.t.\  $U$ is a connected simply--connected
coordinate neighborhood of
a fixed point $(p_0,q_0) \in U$. 
For such a pointed coordinate neighborhood $(U, (p_0,q_0))$ 
with $U \subset U_{(p_0,q_0)}$ we 
define $\tau^U: U \to T_{M^{(1)}}\boxtimes T_{M^{(2)}}|_{(p_0,q_0)}$ by
\begin{equation}
\tau^U(p,q) = (\theta_{\tau })_{(p,q)}^{(p_0,q_0)}(0).
\end{equation}

From now on we will always use pointed neighborhoods and sometimes
drop the explicit mention of the base--point.
\sm

Denote the 
matrix in the chosen basis
defined by $\theta_{\tau}^f$ on an open pointed coordinate neighborhood
$U \in {\Cal W}$
by $\theta_{\tau}^U$. These matrices satisfy
\begin{equation}
(A^{(1)}\otimes A^{(2)})^V_U\theta_{\tau}^U = \theta_{\tau}^V 
(A^{(1)}\oplus A^{(2)})_V^U.
\label{trafo}
\end{equation}
In this notation:
\begin{equation}
\tau^U(p,q) = \theta_{\tau}^U(p,q)
\label{matrix}
\end{equation}
where we identified the point $(p,q)$ with its coordinate vector
in $T_{(p_0,q_0)}$

\sm

Now choose an open pointed cover of ${\Cal U}$ the image of $s$ inside
$M^{(1)}\times M^{(2)}$
subordinate to ${\cal W}$ along which is small enough
for our purposes.
I.e.\ all open sets of the cover are connected simply--connected 
coordinate neighborhoods, their union is a tubular neighborhood of 
$s(\bar M)$ and all intersections of these opens are
connected simply--connected as well. Furthermore,
$\forall (U,(p,q)) \in {\cal U}:  (p,q) \in {\rm Im}(s), U \subset U_{(p,q)}$,
and $\tau^U(U)\in V_{p_iq_i}$.

\sm

For each $U \in {\Cal U}$ we choose $\bar U$ to be a tubular neighborhood
of the image of $U$ inside $V_{pq}$. Again, different choices lead to 
equivalent diagrams.
 
\sm

We obtain the desired tensor product manifold by gluing the open
sets $\bar U$ together using the tensor product of the transition 
functions. More precisely:
\sm
For a given pair $U,V \in {\cal U}, U\cap V \neq \emptyset$ 
we define the affine transformation
\begin{equation}
\bar \varphi_{\bar U}^{\bar V}:= \left( \begin{array}{cc}
\bar A_{\bar U}^{\bar V}& \bar\xi_{\bar U}^{\bar V}\\
0&1
\end{array}
\right)
=
\left( \begin{array}{cc}
A^{(1) V}_U\otimes A^{(2) V}_U& \theta_{\tau}^V (\xi^V_U)\\
0&1
\end{array}
\right)
\label{trans}
\end{equation}
where $A_U^V= A_{(U,(p,q))}^{(V,(p',q'))}$. 
It is straightforward to check using (\ref{equi}) and (\ref{trafo}) that 
\begin{equation}
\bar\varphi_{\bar U}^{\bar V}=\bar\varphi^{\bar U -1}_{\bar V},\quad
\bar\varphi_{\bar V}^{\bar W}\bar\varphi_{\bar U}^{\bar V}=
\bar\varphi_{\bar U}^{\bar W}
\end{equation}
and therefore we get a manifold 
\begin{equation}
N:= (\amalg_{U \in {\Cal U}}\, {\bar U})/{\cal R}
\end{equation}
where ${\cal R}$ is the equivalence relation induced by the $\bar \varphi$.

\sm

The cover $\bar{\cal U}:= \{\bar U|U\in {\cal U}\}$ together
with the inclusion maps $i^{\bar U}:\bar U \to {\Bbb C}^{n_1 n_2}$ and the 
transition
functions $\varphi_{\bar U}^{\bar V}$  (\ref{trans}) define an affine atlas 
of $N$. 
The canonical flat structure $\bigoplus_{a,b} {\Bbb C} \; (\del_a \otimes \del_b)$ 
on each $\bar U$ together with
the tensor metric $g^{\bar U}= g_p \otimes g_q$ glue together under the
affine transformations to form an affine flat structure on $N$.

\sm

Due to the condition (\ref{trafo}) 
the maps $\tau^U:U \to {\bar U} \subset V_{pq}$ of (\ref{matrix}) satisfy:
\begin{equation}
\bar \varphi^{\bar V}_{\bar U}\tau^U(u) = \tau^V(u)
\end{equation} 
for $u \in U \cap V$ and thus glue together to 
a map $t:\bigcup_{U\in{\cal U}} U \to N$. We can now define the map $\tau$ as
\begin{equation}
\tau := t\circ s \circ p
\end{equation}
\sm
Note that since $T_{M^{(1)}}$ and $T_{M^{(2)}}$ are trivial along the
factors ${\Bbb C}$ of the decomposition $M^{(i)}=\bar M^{(i)} \times {\Bbb C}$,
 $s\circ p$ induces the following identity between bundles
\begin{equation}
(s \circ p)^*(T_{M^{(1)}}\boxtimes T_{M^{(2)}}) \cong 
T_{M^{(1)}}\boxtimes T_{M^{(2)}}
\label{piso}
\end{equation}
Due to the fact that all higher (i.e.\ higher than the third) 
derivatives of the potentials
$\Phi^{(i)}$ are independent of the coordinates of the identities,
the above identity of bundles induces an isomorphism of 
germs of pointed Frobenius manifolds
\begin{equation}
(T_{M^{(1)}}\boxtimes T_{M^{(2)}}|_m,0) \cong 
(T_{M^{(1)}}\boxtimes T_{M^{(2)}}|_{s\circ p(m)},0).
\label{giso}
\end{equation}

\sm
So far
we have constructed an affine flat manifold
$N$ a map $\tau:M^{(1)}\times M^{(2)} \to N$ and it is easy to confirm that
$\tau^*(\theta_{can})= \Theta \circ \theta_{\tau}$
and that $\tau^*(T_{N})$ and $T_{M^{(1)}}\boxtimes T_{M^{(2)}}$ 
are naturally isomorphic under the composition $\Theta$ of 
the identification of the tangent space of a vector space at
a point with the vector space itself and the isomorphism of (\ref{piso}).
\sm

Furthermore, the by construction the map $\tau$ factors through $p$ and
since $t$ is an injection on ${\rm Im}(s)$, $i:= t\circ s$ is an embedding.

\sm
To endow $N$ with the desired Frobenius manifold structure, we have
to check that the $\bar U$ glue together as Frobenius manifolds.

This follows from the properties of the tensor product connection.
Let $u\in U \cap V$, $U$ have the base--point $(p,q)$ and
$(p',q')$ be the base--point of $V$.
The Frobenius structure of $\bar U$ is given by the potential $\Phi_{pq}$
and the Frobenius structure of $\bar V$ by the potential $\Phi_{p'q'}$.
Using the tensor connection we see that the germ of 
$\Phi_{pq}$ at $\tau^U(u)$  and
the germ of $\Phi_u$ at $0$ in $V_u$ can be given by the same power series.
This equality between the germs 
holds as well for the germ of $\Phi_u$ at $0$ in $V_u$
and the germ of $\Phi_{p'q'}$ at $\tau^V(u)$
and since the whole germs coincide so do the functions: 
\sm
\begin{equation}
\Phi_{pq}(x)= \Phi_{p'q'}(y) \quad \forall y = 
\bar \varphi^{\bar V}_{\bar U}(x).
\end{equation}

Since $\theta_{\tau}$ is a tensor connection the definition of $t$ and
the isomorphism (\ref{piso}) show that the condition i) for a tensor
product diagram also holds. 

\sm

\pagebreak
{\bf Uniqueness}

\sm

Let two diagrams as in \ref{isodef} be given. 
To construct the open subsets $U_N$,
$U_{N'}$ and the isomorphism
of Frobenius manifolds $\phi$ we choose a pointed cover ${\cal U}_N$ of 
a tubular neighborhood of 
the image of $\tau$ which additionally has the following properties:
\begin{itemize}
\item[i)] $\forall (U,n)\in {\Cal U}_N: n \in {\rm Im}(\tau), U \subset U_n$
where $U_n$ is again the notation for the domain of convergence
of $\Phi^N_n$.
\item [ii)]
The opens and their intersections should be connected and simply--connected and
all open sets should be coordinate neighborhoods.
\item[iii)] Using the isomorphisms $\Theta$ and $\Theta'$ 
we can identify a small neighborhood of 
$0$ on the fiber of $T_{M^{(1)}}\boxtimes T_{M^{(2)}}$ at $(p,q)$
with a neighborhood of $\tau((p,q))$ on $N$ and a neighborhood 
of $\tau'((p,q))$ on $N'$. The condition for our cover is that
all open neighborhoods of the cover are so small that
 the above identifications exist.
\end{itemize}
We define the map $\phi$ as the concatenation of these identifications, i.e.\
\begin{equation}
\phi^{(U,n)}= \exp|_{\tau'(s\circ i^{-1}(n))}
\circ \hat\tau'\circ \Theta'\circ\Theta^{-1}
\circ \hat\tau|_{s\circ i^{-1}(n)}^{-1}\circ\exp^{-1}|_n
\end{equation}
where $\hat\tau|_{s\circ i^{-1}(n)}^{-1}$ is the inverse of
$\hat\tau$ restricted to the fiber of $T_{M^{(1)}\times M^{(2)}}$ 
at $s\circ i^{-1}(n)$.
We set $U_N=\bigcup_{(U,n)\in {\cal U}_N} U$ 
and $U_{N'}= \bigcup_{(U,n)\in {\cal U}_N} \phi^{(U,n)}(U)$.
Since $\tau^*\theta_{can}= \Theta\circ \theta_{\tau}$ and
$\Theta'\circ \theta_{\tau}= \tau^{\prime *}\theta'_{can}$ and all
maps preserve the relevant germs of Frobenius manifolds, 
it is clear that the maps $\phi^{(U,n)}$ patch together as 
a morphism of Frobenius manifolds. On the image of $i$ we have
$\phi|_{{\rm Im}(i)}= \tau'\circ s \circ i^{-1}$. 
A short calculation shows that this map yields a bijection between
${\rm Im}(i)$ and ${\rm Im}(i')$ with 
inverse $\tau\circ s' \circ i^{\prime -1}$. After making
the original cover smaller if necessary, we can assume that the induced cover
${\cal U}_{N'}:= \{(\phi(U), \phi(n) )| (U,n)\in {\cal U}_N\}$ also
satisfies the conditions i)--iii) and the union of the opens of this
cover is again a tubular neighborhood of ${\rm Im}(\tau')$. Hence, 
we can perform the analogous
construction starting from ${\cal U}_{N'}$ yielding 
an inverse morphism. Thus $\phi$ is
an isomorphism of Frobenius manifolds which satisfies
the condition $\phi\circ \tau=\tau'$ and
the commutativity of the upper part of the diagram (\ref{iso}) 
follows directly from the 
construction. 
\sm

\subsubsection{Proposition}
{\it
Given two 
Frobenius manifolds with factorizable flat identities
and Euler fields (with $d=1$) then any tensor product manifold carries
the natural tensor product identity and can be endowed
with an Euler field locally defined
by \ref{ts3.2.1}.}

\sm
{\bf Proof.}
 
We will define  the identity and the Euler field on a tensor on the
product manifold $N$ constructed above. The results can be pushed to any
equivalent manifold. Let
\begin{equation}
e_N^{\bar U} := \del^{(1)}_0\otimes \del^{(2)}_0
\quad E^{\bar U} :=  E_{pq}
\end{equation}
where $E_{pq}$ is the Euler field constructed 
for the formal tensor product in 
\ref{ts3.2.1}. 
The gluing condition for the identities is clear from the
compatibility of the flat structures. What still remains to be shown
is that the locally defined Euler fields glue together
$J_{\bar\varphi_{\bar U}^{\bar V}}E^{\bar U}(x) = 
E^{\bar V}(\bar\varphi_{\bar U}^{\bar V}(x))$, where $J$ is the Jacobian.

Let $\bar\varphi_{\bar U}^{\bar V}(x) = A^{(1)}\otimes A^{(2)} x
+ \theta_{\tau}^V(\xi^V_U)$. Furthermore, we will write
all Euler fields in matrix form; for $U = U^{(1)} \times U^{(2)}$,
$V = V^{(1)} \times V^{(2)}$
$$
E^{U^{(i)}}(x) = D^{U^{(i)}}x + r^{U^{(i)}}, \quad i = 1, 2
$$ 
The gluing conditions for $E^{(1)}$ and $E^{(2)}$ read
$$
A^{(i)} D^{U^{(i)}} = D^{V^{(i)}} A^{(i)}, \quad 
A^{(i)} r^{U^{(i)}} = D^{V^{(i)}}\xi_{U^{(i)}}^{V^{(i)}} + r^{V^{(i)}}, 
\quad i = 1, 2.
$$
In this notation for all $(p,q)$:
$$
E_{pq} = (D^{U^{(1)}} \otimes Id + Id \otimes D^{U^{(2)}} - Id \otimes Id)(x) 
+ \theta_{\tau}^U(r^{U^{(1)}}+r^{U^{(2)}})
$$
And hence
$$
\begin{aligned}
A^{(1)} \otimes A^{(2)} E_{pq}(x) =&
(D^{V^{(1)}} \otimes Id + Id \otimes D^{V^{(2)}} - Id \otimes Id)
(\bar\varphi_{\bar U}^{\bar V}(x))\\
&+ \theta_{\tau}^V \big(-(D^{U^{(1)}}\oplus Id + Id \oplus D^{U^{(2)}}- Id \oplus Id)
(\xi^U_V)+ A^{(1)}r^{U^{(1)}}+A^{(2)}r^{U^{(2)}}\big)\\
=& E_{p'q'}(\bar\varphi_{\bar U}^{\bar V}(x))
\end{aligned}
$$

\sm

\subsubsection{Remarks}
\begin{itemize}
\item[i)]
The restriction of factorizable identity is not too severe. In all
presently known examples this is the case. This includes all
semi--simple Frobenius manifolds considered on ${\Bbb C}^n$ without
the diagonals, as well as the split semi--simple Frobenius manifolds
on the universal cover of the previous space given by special initial
conditions. 
\item[ii)]
Locally one can always complete the direction of the identity by
using an appropriate embedding into ${\Bbb C}^n$ (see below).
\end{itemize}

\subsection{Embedded Frobenius manifolds}

Notice that the universal cover of every affine manifold of dimension $n$
has an immersion into ${\Bbb C}^n$ (cf.\ e.g.\ [KW]). This immersion can
be quite non--trivial however, see e.g.\ [ST]. 
We will call a Frobenius 
manifold an {\it embedded Frobenius manifold} if the manifold itself has
an embedding into ${\Bbb C}^n$. In this case we will identify
the Frobenius manifold with its image under the embedding.
It then has global coordinates given by a choice of basis for the
affine flat tangent bundle.

Actually, most constructions of Frobenius manifolds use global
coordinates e.g.\ the ones coming from quantum cohomology, unfolding
of singularities or
Landau--Ginzburg models.

\sm

\subsubsection{Lemma}
\label{decomp}
{\it An embedded Frobenius manifold ${\bf M}$,
$M \subset {\Bbb C}^n$, with flat identity 
can be completed in the direction
of the identity. I.e.\ there exists a Frobenius manifold $\tilde {\bf M}$ 
with factorizable
flat identity which contains ${\bf M}$ as a Frobenius manifolds
$\tilde M = \bar M \times {\Bbb C} \supset M$ and the other structures
are given by restriction.
}

\sm

{\bf Proof.}
Since the potential a polynomial of order less or equal three in the
coordinate of the identity for a Frobenius manifold with
flat identity, the respective three--tensor defining the Frobenius structure
is independent of the coordinate of the flat identity.
Hence, since the tangent bundle is trivial, 
we can enlarge the domain of definition of this
three--tensor so that it contains all lines in the
direction of the identity. Likewise,
we can extend the metric to these points too.

\sm

The structure of a tensor product manifold for two 
embedded Frobenius manifolds with factorizable flat identity
allows a more explicit description.

\sm

Let $M^{(1)}$ and $M^{(2)}$ be realized in ${\Bbb C}^{n_1}$ respectively ${\Bbb C}^{n_2}$
and consider the map 
$\tau: {\Bbb C}^{n_1+n_2}\to {\Bbb C}^{n_1 n_2}= {\Bbb C}^{n_1}\otimes {\Bbb C}^{n_2}$ given by the matrix
$\tau_{ab}=\del_{0a} + \del_{b0}$. 
Looking at the construction we arrive immediately at the following:

\sm

\subsubsection{Proposition}
{\it 
The tensor product manifold for
embedded Frobenius manifolds is equivalent to an embedded Frobenius manifold
given by a neighborhood
of the image $\tau(M^{(1)}\times M^{(2)}) \subset {\Bbb C}^{n_1}\otimes {\Bbb C}^{n_2}$
and the image is isomorphic to $\bar M^{(1)}\times \bar M^{(2)} \times {\Bbb C}$
in the notation of Lemma \ref{decomp}.
}

\sm

\subsubsection{Remark}
The tensor product in the particular cases considered in
this subsection thus contains a subset isomorphic to the
image of $p$ parameterized by the  
coordinates $x_{a0},x_{0b}$ ($x_{ab}= 0$ for $ab\neq 0$)
The third derivatives of the potential $\Phi$ of the tensor
product satisfies:
\begin{equation}
\Phi_{aa'bb'cc'}|_{{\rm Im}(\tau)}(x_{i0},x_{0j}) = 
\Phi^{(1)}_{abc}(x^{(1)}_i \equiv x_{i0} )
\Phi^{(2)}_{a'b'c'}(x^{(2)}_j \equiv x_{0j}), \quad i,j\neq 0
\label{improd}
\end{equation}
but moreover along the image of $\tau$ the {\em whole germs} of $\Phi$
are given by the tensor product of the associated pointed germs.
Vice versa, the condition (\ref{improd}) does not suffice to 
identify a tensor product manifold, since it does not determine the
higher derivatives in the $\del_{ab}$ directions for $ab\neq 0$.

\subsection{General tensor products}

The construction and the universality statement
easily generalize to the following setting.

\subsubsection{Definition}
Let  $\theta_{\tau}$ be a tensor product connection on the exterior
product bundle over a Cartesian product of two 
Frobenius manifolds  ${\bf M}^{(1)}$ and ${\bf M}^{(2)}$.
We call the Cartesian product ${\bf M}^{(1)}\times {\bf M}^{(2)}$
{\em $\theta_{\tau}$--reducible} if there exists a flat affine manifold
$\bar M$, an affine projection 
$p:M^{(1)}\times M^{(2)}\to \bar M$ together with an affine section $s$ of this
projection which is an embedding of affine flat manifolds 
satisfying the following conditions
\begin{itemize}
\item[i)] The projection condition of (\ref{piso})
\begin{equation}
(s \circ p)^*(T_{M^{(1)}}\boxtimes T_{M^{(2)}}) \cong
T_{M^{(1)}}\boxtimes T_{M^{(2)}}
\end{equation}
\item[ii)] The condition on the respective germs of Frobenius manifolds
(\ref{giso})
\begin{equation}
(T_{M^{(1)}}\boxtimes T_{M^{(2)}}|_{s\circ p(m)},0)
\stackrel{(s\circ p)^*}{\cong}
(T_{M^{(1)}}\boxtimes T_{M^{(2)}}|_m,0) 
\end{equation}
\item[iii)] The compatibility with the tensor product connection;
under the identification of i):
\begin{equation}
(s \circ p)^*\theta_{\tau} = \theta_{\tau}
\end{equation}

\item[vi)] And the embedding condition
\begin{equation}
T_{M^{(1)}} \boxplus T_{M^{(2)}}|_{{\rm Im}(s)} = 
T_{{\rm Im}(s)} \oplus ker(\theta_{\tau}^f)|_{{\rm Im}(s)}
\end{equation}
\end{itemize}

The triple $(\bar M,p,s)$ is then called {\em a $\theta_{\tau}$--reduction}

\sm
\subsubsection{Lemma}
{\it
If in a Cartesian product $M^{(1)}\times M^{(2)}$ each of the factors $M^{(i)}$
can be decomposed as $\bar M^{(i)}\times {\Bbb C}^{n_i}$ as flat
affine manifolds and the whole Cartesian product can be decomposed 
as $\bar M\times {\Bbb C}^n$ 
where the third derivatives of the potential are constant in 
the ${\Bbb C}$--directions and the kernel of $\theta_{\tau}$ gives
the coordinates the ${\Bbb C}$--directions and furthermore the factor
${\Bbb C}^n$ is an affine flat factor 
of ${\Bbb C}^{n_1}\times {\Bbb C}^{n_2}$,
then ${\bf M}^{(1)}\times {\bf M}^{(2)}$ is
$\theta_{\tau}$--reducible.}

\sm

{\bf Proof.}
We can decompose $M^{(1)}\times M^{(2)}$ as $\bar M^{(1)}\times \bar M^{(2)}
\times {\Bbb C}^{n_1+n_2-n} \times  {\Bbb C}^n$ where 
${\Bbb C}^{n_1+n_2-n} \times {\Bbb C}^n$ is the postulated decomposition
of ${\Bbb C}^{n_1+n_2}$.
Now consider the projection 
$p: M^{(1)}\times M^{(2)} \to \bar M= \bar M^{(1)}\times \bar M^{(2)}
\times {\Bbb C}^{n_1+n_2-n}$
and the zero section $s$. They satisfy all the conditions of 
a $\theta_{\tau}$--reduction.
\sm

\subsubsection{Lemma}
\label{locgen}
For two Frobenius manifolds with not necessarily factorizable flat identities
a general tensor product diagram for the canonical tensor product connection
${\theta}_{\tau}$ (\ref{thetatau}) exists.

\sm
{\bf Proof.}
Locally we can achieve the following situation:
let  $U$ in $M^{(1)} \times M^{(2)}$ be an open set which satisfies:
$U= U^{(1)}\times U^{(2)} = \bar U^{(1)}\times D_r \times \bar U^{(2)} 
\times D_r$
where $D_r$ is a disc of radius $r$ in ${\Bbb C}$ centered at $0$ 
---$D_r=\{z \bigm| |z| < r\}$--- which
is coordinatized by the identity.
Now take $\bar M := \bar U^{(1)}\times \bar U^{(2)} \times D_{2r}$ where
$D_{2r} = \{z\bigm| |z| < 2r \}$ and define $p$ and $s$ as in
the construction of \ref{exist}. Since again the potential is constant 
on the factors $D$ 
all necessary properties directly follow.

\sm

\subsubsection{Definition\label{gendef}}
For three Frobenius manifolds 
${\bf M}^{(1)}, {\bf M}^{(2)}$ and ${\bf N}$ consider a diagram of the type
\begin{equation}
\xymatrix{
T_{M^{(1)}} \boxtimes T_{M^{(2)}}\ar[r]^-{\Theta}\ar[d]&
\tau^*T_N\ar[dr]\ar[r]^{\hat\tau}&T_N\ar[d]\\
M^{(1)}\times M^{(2)}\ar[d]^p\ar[rr]^{\tau}&&N\\
\bar M\ar@/^1pc/[u]^{s}\ar@{^{(}->}[urr]^{i}&&\\
}
\end{equation}
together with a tensor product connection $\theta_{\tau}$, such that
${\bf M}^{(1)}\times {\bf M}^{(2)}$ is $\theta_{\tau}$--reducible, 
$(\bar M,p,s)$ is a $\theta_{\tau}$--reduction, $\tau$ is an affine map
and $\Theta$ is an isomorphism of metric bundles
with affine flat structure.

We will call such a diagram a {\em $\theta_{\tau}$--tensor product diagram} 
if it satisfies the condition i) and ii) of a tensor product diagram where
now $\theta_{\tau}$ is given in the data.
\sm

\subsubsection{Definition}
Two general tensor product diagrams are called {\em equivalent}
if there exist open neighborhoods $U_N, U_{N'}$ of the images of 
$\tau, \tau'$
and an isomorphism $\phi$ of
of Frobenius manifolds $U_N\to U_{N'}$ s.t.\ the induced diagram of the
form (\ref{iso})
satisfying the conditions
of a general tensor product diagram is commutative.

\subsubsection{Theorem}
{\it

Given a tensor product connection $\theta_{\tau}$ on an exterior
product bundle over the Cartesian product of two  
Frobenius manifolds  ${\bf M}^{(1)}$ and ${\bf M}^{(2)}$ 
with $\theta_{\tau}$--reducible
Cartesian product,
there exists a 
$\theta_{\tau}$--tensor product diagram and thus a tensor product manifold.
Furthermore, fixing a $\theta_{\tau}$--reduction all diagrams involving
this reduction are equivalent. 
}

\sm

{\bf Proof.}
In the case that and
$\theta_{\tau}$ 
is not identically zero we can retrace the proof of 
\ref{exist}, 
since we only used that
$\theta_{\tau}$ is a non--zero tensor product connection and the existence
of a $\theta_{\tau}$--reduction. 
In case $\theta_{\tau} \equiv 0$ the image of $p$ is just
a point and the tensor product already exists by the construction
for convergent germs of Frobenius manifolds. The uniqueness then
follows directly from the condition i) for a general tensor product diagram.

\sm

\subsubsection{Remark} 
The local situation for Frobenius manifolds with flat identities can be
described in three different ways:
\begin{itemize}
\item[i)] Via the {\em tensor product connection}.\\
Given any point $(p,q)\in M^{(1)} \times M^{(2)}$ we have the corresponding
 pointed germ of the tensor product relative to the pair of 
base--points $(p,q)$ on the fiber 
of $T_{M^{(1)}} \boxtimes T_{M^{(2)}}$
over $(p,q)$. This germ can be seen as the germ describing the situation locally
since by virtue of the existence of the tensor product connection 
of \ref{thetatau} all continuations of this germ will contain all neighboring
germs.
\item[ii)] Via an {\em embedding and completion}.\\
Locally we can embed any complex 
affine manifold of dimension $n$ into ${\Bbb C}^n$.
Using the analysis of embedded Frobenius manifolds with flat identity
we can complete our embedded manifold and use the Theorem for tensor product
diagrams to find a local tensor product manifold.
\item[iii)] Via the {\em local general tensor product}.\\
See Lemma \ref{locgen}.
\end{itemize}
Of course all these descriptions are compatible. The compatibility of i)
with ii) and iii) is manifest in the condition i) of the 
definition of (general)
tensor product diagrams. For a suitable neighborhood, 
we can pass from the ii) to iii)  by restricting everything 
to the image of this neighborhood in the completion of its embedding
and its image under $\tau$.

\sm

\subsubsection{Remark}
One necessary condition for the existence of 
a general tensor product diagram is the existence
of a tensor product connection. In some cases there are many such
connections in others there may be none. If there is none one this shows
that the tensor product can only defined w.r.t.\ a fixed base--point
and there is no way of naturally parameterizing the
tensor products by a submanifold in a Frobenius manifold. 

\sm

\subsubsection{Examples}
\begin{itemize}
\item[i)]
Choose two vector spaces $V,V'$ with a constant product and consider the
constant tensor product multiplication in the sense of algebras 
on $V\otimes V'$. Now any linear map 
$\theta_{\tau}:V\times V' \to V\otimes V'$ will provide a tensor product
connection.
\item[ii)]
Consider two one--dimensional Frobenius manifolds whose potentials
both have a zero; say at the points $p$ and $q$, but are not
constantly zero near these points. Using the Proposition \ref{zerotensor}
we see that the Frobenius structure
on the tangent space at $(p,q)$ would be 
given by a vanishing
potential. On the other hand, near the point $(p,q)$ the 
potentials are non--vanishing by assumption and so is their product
which is the value of $\Psi$ on the zero section of the exterior product
bundle. Thus there is no tensor
product connection on a  
neighborhood of the point $(p,q)$. (For a general statement about
the one--dimensional situation see below.)
\item[iii)]
A tensor product connection
always exists for the product of any Frobenius manifold 
${\bf M}$ with  ${\Bbb C}$ carrying a constant multiplication, i.e.\ 
the third derivative potential $\Phi$ is constant 
$\Phi_{zzz} \equiv \alpha$ ---
where $z$ is a fixed coordinate on ${\Bbb C}$.
Notice that after scaling $\frac{\del}{\del z}$ we can assume
that the multiplication is either constantly zero
or $\frac{\del}{\del z}$ is a flat identity for ${\Bbb C}$.
In this case, $\theta_{\tau}(\del_a)= \del_{a0}, \theta_{\tau}(\del'_0) = 0$
provides a tensor product connection. 
The tensor product manifold is $M$ itself and the
map $\tau=\pi_1$  the first 
projection of $M\times {\Bbb C}$. Here $p=\tau$ and $s$ is the zero
section. 
\end{itemize}
\sm

\subsubsection{The tensor product of one--dimensional Frobenius manifolds}

In this subsection we give a complete answer to the existence question of
a tensor product connection
in the case of two one--dimensional Frobenius manifolds.

\med

\subsubsection{Proposition}
{\it
A tensor product connection for two one--dimensional Frobenius 
manifolds exists
if and only if one of the factors has a locally constant multiplication i.e.\
locally $\Phi_{zzz} = \alpha$ where $z$ is the coordinate function of 
the flat vector field 
$\frac{d}{dz}$
or equivalently it carries a flat identity or a zero multiplication.
} 

In other words, the tensor product of 
one--dimensional theories is essentially pointed and
does not contain perturbations of the base--points.
\sm

{\bf Proof.}
Using the notation of
the previous subsection, suppose a tensor product connection exists.
Choose any point $(p,q)\in M^{(1)} \times M^{(2)}$ and a coordinate
neighborhood $U$ of $(p,q)$ with the local normalized
product coordinate $(z_1,z_2)$, i.e.\ $(z_1,z_2)(p,q)=(0,0)$.
Write $\varphi_1 := \Phi^{(1)}_{z_1z_1z_1}$
and  $\varphi_2 := \Phi^{(2)}_{z_2z_2z_2}$ and
$\varphi_{z_1,z_2}:=\frac{\del^3}{\del^3z}\Phi_{z_1} \otimes \Phi_{z_2}$

Set
$$
\psi(z_1,z_2,z) := \varphi_{z_1,z_2}(z- \tau(z_1,z_2)).
$$
The function $\psi$ is independent of $z_1,z_2$ since by definition of
a tensor product connection
$$
\psi(z_1,z_2,z) = \varphi_{0,0}(z)
$$
Therefore:
$$
\frac{\del}{\del z_i} \psi \equiv 0\\
$$
and thus 
\begin{gather*}
\frac{\del}{\del z_i} \psi(z_1,z_2,\tau(z_1,z_2))=
\frac{\del \varphi_{z_1,z_2}}{\del z_i}(0)
-\frac{\del \tau}{\del z_i}|_{(z_1,z_2)} \varphi_{z1,z2}'(0)=0\\
\Leftrightarrow  \frac{\del}{\del z_i}(\varphi_1 \varphi_2)- 
\frac{\del \tau}{\del z_i}|_{(z_1,z_2)} 
(\varphi_1' \varphi_2^2 + \varphi_1^2 \varphi_2')=0
\end{gather*}
where $'$ denotes the derivative and we used the short--hand notation
$\varphi_1,\varphi_2$ for $\varphi_1(z_1),\varphi_2(z_2)$.
Furthermore,
\begin{equation}
\begin{gathered}
\frac{\del}{\del z_i}\psi'(z_1,z_2,\tau(z_1,z_2))=
\frac{\del \varphi'_{z_1,z_2}}{\del z_i}(0)
-\frac{\del \tau}{\del z_i}|_{(z_1,z_2)} \varphi_{z1,z2}''(0)=0\\
\Leftrightarrow  \frac{\del}{\del z_i}(\varphi_1' \varphi_2^2 + \varphi_1^2 
\varphi_2')
- \frac{\del \tau}{\del z_i}|_{(z_1,z_2)} 
(\varphi_1'' \varphi_2^3 + 5\varphi_1 \varphi_2\varphi_1' \varphi_2'
+\varphi_1^3\varphi_2'' )=0.
\label{cond2}
\end{gathered}
\end{equation}

Therefore we have
$$
((\varphi_2\frac{\del}{\del z_1}+\varphi_1\frac{\del}{\del z_2})
\psi'(z_1,z_2,\tau(z_1,z_2)))\psi'(\tau(z_1,z_2))
= \varphi_1 \varphi_2\varphi_1' \varphi_2'
(\varphi_1' \varphi_2^2 + \varphi_1^2 \varphi_2')=0.	
$$

Therefore either $\varphi_1$ or $\varphi_2$ constantly vanish, or  
we may assume
that on some open set $\varphi_1(z_1) \neq 0$ and $\varphi_2(z_2) \neq 0$
and therefore if also neither $\varphi'_1$ nor $\varphi'_2$ constantly vanish
we must have
\begin{equation}
\frac{\varphi_1'(z_1)}{\varphi_1(z_1)} = 
- \frac{\varphi_2'(z_2)}{\varphi_2(z_2)}= c
\label{diffeq}
\end{equation}
where $c$ is a constant.

\sm

The solution to these simple differential 
equations is $\varphi_1(z_1)= \frac{1}{cx+d_1}$
and $\varphi_2(z_2)= \frac{2}{cx+d_2}$. 
Inserting (\ref{diffeq})
into the equation (\ref{cond2}) we find
$$
c^3 \varphi_1^5 \varphi_2^4 \equiv 0
$$
which yields that $c=0$ a contradiction to the last assumption.

\sm
Therefore either  $\varphi'_1\equiv 0$ or $\varphi'_2 \equiv 0$ and the
proposition follows.
\section{\label{ts5.1}Semi--simple Frobenius manifolds} 

\sm
\subsection{Semi--simple Frobenius manifolds} 

\sm

We will briefly recall the main notions of semi--simple Frobenius
manifolds as
explained in [M1]. For other versions see [D2] or [H].
A Frobenius manifold of dimension $n$ is called semi--simple
(respectively split semi--simple), 
if an isomorphism
of the sheaves of ${\Cal O}_M$--algebras
\begin{equation}
({\Cal T}_M, \circ)  \simeq ({\Cal O}_M^n,\text{ componentwise
multiplication})
\label{t5.1}
\end{equation}
exists everywhere locally (respectively globally).

\sm
If a Frobenius manifold $M$ is semi--simple, one can find so--called
canonical coordinates 
$u_i$ ---unique up to constant shifts and renumbering---
 s.\ t.\ the metric and the three--tensor $A$ defining the
multiplication become
particularly simple. Let $e_i =  \frac{\del}{\del u_i}, \nu_i =
{\mathrm d}u_i$, then
\begin{equation}
g= \sum_i \eta_i (\nu_i^2), \quad
A=\sum_i \eta_i(\nu_i^3).\label{t5.2}
\end{equation}

\sm

If in addition an Euler field exists, then it has the form
$E = \sum (u^i +c^i) e_i.$
In this situation, we will normalize the coordinates in such a way
that
\begin{equation}
E = \sum u^i e_i.
\label{t5.3}
\end{equation}
This normalization fixes the ambiguity in the coordinates $u^i$ and
renders them unique up to the ${\Bbb S}_n$--action.

\sm

\subsubsection{Definition}
In the above situation, we will call a point $m \in M$ {\em tame}, if it 
satisfies $u_i(m)\neq u_j(m)$ for all $i\neq j$. 
In other words, the point $m$ is tame, if the spectrum of 
the operator $E\circ$ on ${\Cal T}_M$ is simple. 

\sm

In the theory of semi--simple Frobenius manifolds one then defines 
certain natural structure connections which give rise to isomonodromic
deformations.
These deformations 
are governed by the Schlesinger differential 
equations [Sch, Mal], thus providing
a link between Frobenius manifolds and solutions of the Schlesinger
equations; the details can be
found in [D2, M1, MM].

\sm

\subsubsection{\label{ts5.1.6}Theorem (2.6.1 of [MM])}
{\it
Let $(M,(u^i),T,(A_i))$ be a strictly special solution and $e$ an
identity of weight $D$,
then these data come from a unique structure of semi--simple split
Frobenius manifold $M$
with an identity $(d_0 =1)$ and an Euler field via 
\begin{equation}
\begin{gathered}
T = \Gamma(M, {\Cal T}_M^f),\quad
(u^i)\text{: the canonical coordinates}\\
A_j (e_i) =0\quad \text { for $ i\neq j$},\quad
A_i(e_i) = - \frac{1}{2} e_i + \sum_{j:j \neq i} (u^j - u^i)
\frac{\eta_{ij}}{\eta_i}e_j\\
\text {The operator $\Cal V$ is given by: }
{\Cal V}(X) = \nabla_{0,X}(E) -\frac {D}{2}X. 
\label{t5.11}
\end{gathered}
\end{equation}
Here, the manifold $M$ only has tame points which means that by definition
$u^i (m) \neq u^j(m), \forall i\neq j ,m \in M$. $M$ should be
regarded as a splitting cover of the subspace of tame points of 
a given Frobenius manifold.
}

\sm
For the notion of strictly special solutions consult [MM].
\sm

\subsubsection{\label{ts5.1.7}Special initial conditions}

\sm

Fixing a base--point in a solution to Schle\-singer's equations and 
taking the coordinates $e_i$ for $T$ call
a family of matrices
$A_1^0, \dots A_m^0 \in {\mathrm End}(T)$ {\it special initial
conditions}, if there exists
a diagonal metric $g$ and a skew--symmetric operator ${\Cal V}$
s.\ t.\ 
$A_j^0 = - ({\Cal V} + \frac {1}{2}{\mathrm Id}) P_j$, where $P_j$
is the projector
onto ${\Bbb C} e_j$. 

In the case of semi--simple Frobenius manifolds with
an Euler field and a flat identity,
the special initial conditions are given by the
value of the structures listed in \ref{ts5.1.6} at a fixed tame 
point $m_0 \in M$ with coordinates
$(u^i_0)$; more precisely, the metric is given by 
the $\eta_i (m_0)$ and the operator 
${\Cal V}$ by the matrix $(v_{ij})_{ij}$ defined 
by $(\nabla_{0,e_i}(E) -\frac {D}{2} (e_i)) (m_0)= (\sum_j v_{ij}e_j)(m_0)$.
The matrix coefficients $v_{ij}$ can be calculated as follows: 
\begin{equation}
v_{ij}= (u^i-u^j)\frac{\eta_{ij}}{\eta_i}(m_0).
\end{equation}
 
\sm

\subsection{\label{ts5.2}The tensor product for 
split semi--simple Frobenius
manifolds with Euler field and flat identity}

\sm

In the previous section we constructed
a tensor product of Frobenius manifolds with factorizable flat
identity. 
Moreover any split semi--simple Frobenius manifold is already 
determined by the special initial conditions
at a tame semi--simple point.

\sm

If there is a pair of tame semi--simple 
points $(p,q)\in M^{\prime}\times M^{\prime\prime }$, 
then the image $\tau((p,q))$ is again a semi--simple point, 
since the algebra in the tangent space over the
 base--point of the
tensor product is just the tensor product of two semi--simple
algebras and thus it is itself semi--simple. 

\sm

Thus, the tensor product manifold 
is given locally near $\tau((p,q))$ 
by the special initial conditions at $\tau((p,q))$ if
this new base--point is again tame. The tensor product is globally given by
these conditions for the tensor product of two split semi--simple 
Frobenius manifolds. 

\sm

This condition, however, is
not very restrictive and if there is a pair of tame semi--simple points in 
some open $U$
one can always find a pair of
tame semi--simple points $(p',q')\in U$ whose image is also
tame semi--simple as we will show later.

\med

\subsubsection{\label{ts5.2.3}Canonical coordinates}

\sm

Since the proof of existence of the tensor product makes extensive
use of the flat coordinates, a
natural question to ask  in the setting of semi--simple
Frobenius manifolds is: Is there also a nice formulation
in terms of canonical
coordinates? Generally, one can not expect simple formulas, since
the algebra
in the tangent space over a given point in the tensor product
manifold is generally not
a tensor product of algebras ---this locus is described by the image
of the Cartesian product--- and the ``coupling'' of
algebras results in a destruction of the pure tensor form of
the idempotents.

\med
 
Using the definitions of the tensor product for formal Frobenius manifolds,
we can, however, calculate 
the idempotents of the tensor product in terms of flat coordinates
in the formal situation.
They are given by the following Proposition up to terms of order two in
flat coordinates which is the precision needed to calculate
the special initial conditions.

\med
\subsubsection{\label{ts5.2.4}Proposition}
{\it
Given two semi--simple Frobenius manifolds ${\bold M'},{\bold M''}$ 
let the idempotents near the base--points $m'_0,m''_0$ have the
expansions
$e'_i = e_i^{\prime 0} + \sum x^{\prime a'} e_i^{\prime a'} + O(x^{\prime 2})$,
and 
$e''_i = e_i^{\prime\prime  0} + \sum x^{\prime\prime  a''} 
e_i^{\prime\prime  a''} +
O(x^{\prime\prime  2})$ in the flat coordinates $x'$ and $x''$,
then the idempotents $e_{ij}$ of the tensor product 
$({\bold M'}\otimes_{(m'_0,m''_0)}{\bold M''},0)$ have the
following expansion in the flat coordinates $x$ around $0$:
\begin{equation}
e_{ij}(x) = e_i^{\prime 0}\otimes e^{\prime\prime  0}_j + \sum_{a',a''}
x^{a'a''} (\la_j^{a''} e_i^{\prime a'} \otimes e^{\prime\prime  0}_j
+ \la_i^{a'} e_i^{\prime 0} \otimes  e^{\prime\prime  a''}_j) + O(x^2)
\label{t5.12}
\end{equation}
where $\del'_{a'} = \sum \la_i^{a'} e_i^{\prime 0}$ and 
$\del''_{a''}= \sum \la_j^{a''} e^{\prime\prime  0}_j$.

\sm

Furthermore, the respective coordinate
functions for the tensor
metric $\eta_{ij}:= \eta (e_{ij},e_{ij})$  have the expansions:
\begin{gather}
\eta_{ij}(x)= \eta'_i(m'_0) \eta''_j(m''_0) + \sum_{a',a''} x^{a'a''} 
(\la_j^{a''} (\del'_{a'}\eta'_i)(m'_0)\eta''_j(m''_0) +
\la_i^{a'}\eta'_i(m'_0)(\del''_{a''}\eta''_j(m''_0)))
\label{t5.13}\\
\intertext {and their derivatives 
$\eta_{ij,kl} := e_{kl} \eta_{ij}$ have the following values 
at the base--point $m_0$:}
\eta_{ij,kl}(0) = \delta_{j,l} 
\eta'_{ik}(m'_0)
\eta''_j(m''_0)
+ \delta_{i,k}  \eta'_i(m'_0) 
\del''_{a''} \eta''_{jl}(m''_0)
\label{t5.14}  
\end{gather}
where $\delta_{i,k}$ is the Kronecker delta symbol.

\sm
If the factors carry flat identities and
Euler fields, then the normalized canonical coordinates of $m_0$ are:
\begin{equation}
u^{ij}(0) = u^{\prime i}(m'_0) + u^{\prime\prime  j} (m''_0). \label{t5.15}
\end{equation}
}

\sm

{\bf Proof.} To check the formula (\ref{t5.12}), expand the potential
$\Phi$ up to order
four in the flat coordinates and verify the idempotency by direct calculation.

\sm

The equations for the 
idempotents $e_i= e_i^0 + \sum x^a e_i^a + {\rm O}(x^2)$ 
in flat coordinates are in zeroth order:
\begin{equation}
e_i^0=(e_i^0, e_i^0) 
\label{zeroord}
\end{equation}
and in first order
\begin{equation}
e_i^a= (e_i^0,e_i^0,\del_a) + 2 (e_i^a,e_i^0).
\label{firstord}
\end{equation}
Here we used the notation $(\del_i,\dots,\del_j)$ 
for the higher order multiplications.

\sm

We can now check both conditions for (\ref{t5.12}). For the zeroth order
we obtain:
$$
(e_{ij}^0, e_{ij}^0) =
(e_i^{\prime 0},e_i^{\prime 0}) \otimes
(e_j^{\prime\prime  0},e_j^{\prime\prime  0})=
e_i^{\prime 0}\otimes e_j^{\prime\prime  0}= e_{ij}^0
$$

And for the first order terms:
$$
\begin{aligned}
(e_{ij}^0,& e_{ij}^0,\del_{a'a''}) + 2 (e_{ij}^{a'a''}e_{ij}^0)\\
=&(e_i^{\prime 0},e_i^{\prime 0},\del_{a'})
\otimes ((e_j^{\prime\prime  0},e_j^{\prime\prime  0})\del_{a''})
+((e_i^{\prime 0},e_i^{\prime 0})\del_{a'})
\otimes (e_j^{\prime\prime  0},e_j^{\prime\prime  0},\del_{a''})\\
&+2\lambda_j^{a''}(e_i^{\prime a'},e_i^{\prime 0}) \otimes
(e_j^{\prime\prime  0},e_j^{\prime\prime  0})
+2\lambda_i^{a'}(e_i^{\prime 0},e_i^{\prime 0})
\otimes (e_j^{\prime\prime  a''},e_j^{\prime\prime  0})\\
=& \lambda_j^{a''} 
[(e_i^{\prime 0},e_i^{\prime 0},\del_{a'})+2(e_i^{\prime a'},e_i^{\prime 0})]
\otimes e_j^{\prime\prime  0}
+\lambda_i^{a'} e_i^{\prime 0}\otimes 
[(e_j^{\prime\prime  0},e_j^{\prime\prime  0},\del_{a''})
+2(e_j^{\prime\prime  a''},e_j^{\prime\prime  0})]\\
=&
\lambda_j^{a''} e_i^{\prime a'}\otimes e_j^{\prime\prime  0}+
\lambda_i^{a'} e_i^{\prime 0}\otimes e_j^{\prime\prime  a''}\\
=&e_{ij}^0 
\end{aligned}
$$
since $((e_j^{\prime\prime  0},e_j^{\prime\prime  0})\del_{a''})=
(e_j^{\prime\prime  0},\del_{a''})= \lambda_j^{a''}e_j^{\prime\prime  0}$
and $((e_i^{\prime 0},e_i^{\prime 0})\del_{a'})= (e_i^{\prime 0},\del_{a'})
=\lambda_i^{a'}e_i^{\prime 0}$.

\sm

The expansion for the metrics of the factors reads:
$$
\begin{aligned}
\eta_i'(x') &= \eta'_i(m'_0) + \sum x^{\prime a'} 
2 g(e^{\prime a'}_i,e^{\prime 0}) + O(x^{\prime 2}) \\
\eta_j''(x'') &= \eta''_j(m''_0) + \sum x^{\prime\prime  a''} 
2 g(e^{\prime\prime  a''}_i,e^{\prime\prime  0}) +O(x^{\prime\prime  2})
\end{aligned}
$$

Inserting (\ref{t5.12}) into the tensor metric we obtain (\ref{t5.13}):
\begin{align*}
\eta_{ij}(x) &= g(e_{ij},e_{ij}) \hfill\\ 
&=g(e_i^{\prime 0}\otimes e^{\prime\prime  0}_j,
e_i^{\prime 0}\otimes e^{\prime\prime  0}_j)\\
&\quad + 2 \sum_{a',a''} x^{a'a''} 
g(\la_j^{a''} e_i^{\prime a'} \otimes e^{\prime\prime  0}_j
+ \la_i^{a'} e_i^{\prime 0} \otimes  e^{\prime\prime  a''}_j,
e_i^{\prime 0}\otimes e^{\prime\prime  0}_j)  + O(x^2)\\
&=g'(e_i^{\prime 0}\otimes e_i^{\prime 0})
g''(e^{\prime\prime  0}_j,e^{\prime\prime  0}_j) \\
&\quad + 2 \sum_{a',a''} x^{a'a''}
[\la_j^{a''}g'(e_i^{\prime a'},e_i^{\prime 0})
g''(e^{\prime\prime  0}_j,e^{\prime\prime  0}_j)
+ \la_i^{a'}g'(e_i^{\prime 0},e_i^{\prime 0})
g''(e^{\prime\prime  a''}_j,e^{\prime\prime  0}_j)]
+ O(x^2)\\
&=\eta'_i(m'_0) \eta''_j(m''_0) + \sum_{a',a''} x^{a'a''} 
(\la_j^{a''} (\del'_{a'}\eta'_i)(m'_0)\eta''_j(m''_0) +
\la_i^{a'}\eta'_i(m'_0)(\del''_{a''}\eta''_j)(m''_0))
\end{align*}

The formula (\ref{t5.14})
then follows by derivating (\ref{t5.13}) w.r.t.\
$e_{kl} = \sum_{a'a''} \la_{a'}^i \la_{a''}^j$
where $(\la_{a'}^i)$ and  $(\la_{a''}^j)$ are the inverse matrices
of $(\la^{a'}_i)$ and  $(\la^{a''}_j)$ defined above.
$$
\begin{aligned}
(e_{kl} \eta_{ij})(0) &= (\sum_{a'a''} \la_{a'}^k \la_{a''}^l \del_{a'a''}\eta_{ij})(0)\\
&= \sum_{a'a''} \la_{a'}^k \la_{a''}^l 
(\la_j^{a''} (\del'_{a'}\eta'_i)(m'_0)\eta''_j(m''_0) +
\la_i^{a'}\eta'_i(m'_0)(\del''_{a''}\eta''_j)(m''_0))\\
&= \delta_{jl} 
(\sum_{a'}\la_{a'}^k\del'_{a'}\eta'_i)(m'_0)\eta''_j(m''_0)
+ \delta_{ik}  
\eta'_i(m'_0)(\sum_{a''}\la_{a''}^l \del''_{a''}\eta''_j)(m''_0)\\
&=\delta_{j,l} 
\eta'_{ik}(m'_0)
\eta''_j(m''_0)
+ \delta_{i,k}  \eta'_i(m'_0) 
\del''_{a''} \eta''_{jl}(m''_0)
\end{aligned}
$$

Finally,
(\ref{t5.15}) can be derived
from the expansion of the equation $E = \sum u^{ij}e_{ij}$ with
the Euler field
$E$ given by Theorem \ref{ts4.2.1}. 
$$
\begin{aligned}
\sum_{ij} u^{ij}(0) e_{ij}^0 &=  E(0)
= E'(m'_0)\otimes \del''_0 +\del''_0 \otimes E''(m''_0)\\
&= (\sum_i u^{\prime i}(m'_0) e^{\prime 0}_i) \otimes 
(\sum_j e^{\prime\prime  0}_j)
+ (\sum_i e^{\prime 0}_i) \otimes 
(\sum_j u^{\prime\prime  j}(m''_0) e^{\prime\prime  0}_j)\\
&= \sum_{ij} (u^{\prime i}(m'_0)+  u^{\prime\prime  j}(m''_0))  e_{ij}^0 
\end{aligned}
$$
\med

\subsubsection{Remark}
Notice that since $\del_0 = \sum_i e_i$ we always have
that
$\la_i^0=\la_j^0=1$, so that we retrieve the previous result for the
tensor product constructed in the last section
$$
e_{ij}|_{{\rm Im}\tau} = e_i \otimes e_j
$$
up to order two as it should be. 

\subsection{\label{ts5.2.5}Tensor product of special initial conditions}

\sm

In the Lemma \ref{ts5.2.4},
we have calculated all of the structures (\ref{t5.11})
necessary to determine the special initial conditions.
\sm
We find:
\begin{equation}
\begin{aligned}
v_{ij,kl}&= (u^{ij}-u^{kl})\frac{\eta_{ij,kl}}{\eta_{kl}}\\
&=(u^i+u^j -u^k-u^l)
\frac{\delta_{jl}\eta'_{ik}\eta''_j + \delta_{ik}\eta'_i \eta''_{jl}}
{\eta'_k\eta'_l}\\
&= \delta_{jl}(u^i-u^k)\frac{\eta'_{ik}}{\eta'_k}
+ \delta_{ik}(u^j -u^l)\frac{\eta''_{jl}}{\eta'_l}\\
&=\delta_{jl}v_{ik}+\delta_{ik}v_{jl}
\end{aligned}
\end{equation}
\med

Another approach using the Euler field is 
given by the following observation:

\med

\subsubsection{\label{ts5.2.6}Remark}

\sm

To give the special initial conditions for a tensor product with
the
choice of tensor metric and the Euler field (\ref{t3.7}), 
it suffices to determine the operator
\begin{equation}
{\Cal V}: {\Cal V}(X) = \nabla_{0,X}(E)- \frac{D}{2}X
\label{t5.16}
\end{equation} 
in the tangent space to the base--point ${\Cal T}_{M,m_0}$.
Since ${\Cal V}$ is an ${\Cal O}_{M}$--linear tensor,
its value on a vector field $X$ 
is already determined by $X\bigm|_{m_0} \in {\Cal T}_{M,m_0}$, so that, 
if we are only interested in the operator ${\Cal V}$
restricted to ${\Cal T}_{M,m_0}$, we can use any extension 
of the vector  $X\bigm|_{m_0}$ to a vector field in a neighborhood of $m_0$. 
Choosing a flat extension $X^f$,
the formula (\ref{t5.16}) simplifies to
\begin{equation}
{\Cal V}(X)\bigm|_{m_0} = ([X^f,E] - \frac{D}{2}X^f)\bigm|_{m_0}.
\label{t5.17}
\end{equation}

\sm

In particular, in the situation of Theorem \ref{ts3.2.1}, we can 
extend the
idempotents $e_{ij}\bigm|_{m_0}$ to
flat vector fields $e_{i,j}^f$  and use the formula (\ref{t5.17}) 
to calculate the special initial conditions
via the operator ${\Cal V}$ for the semi--simple tensor Euler field. 
Now it is clear that 
$(e_{ij})\bigm|_{m_0} = e'_i \bigm|_{m'_0} \otimes \; e''_j\bigm|_{m''_0}$,
since the algebra over $m_0$ is just the tensor of the algebras
at the chosen zeros $m'_0$ and $m''_0$. 
Recalling the form of $E$ given by (\ref{t3.7}), we find for
flat $X,Y$
\begin{equation}
[X\otimes Y,E] = [X,E']\otimes Y + X\otimes[Y,E''] - d X\otimes Y.
\label{t5.18}
\end{equation}
Thus,
\begin{multline}
{\Cal V}(e_{ij}^f) = [e_{i,j}^f,E] -\frac {D}{2}e_{i,j}^f =\\
([e_i^{\prime f},E']- \frac {D'}{2}e_i^{\prime f})
\otimes e_j^{\prime\prime  f} +
e_i^{\prime f} \otimes ( [e_j^{\prime\prime  f},E''] 
- \frac {D''}{2} e_j^{\prime\prime  f}).
\label{t5.19}
\end{multline}

\med

\noindent Using the explicit formulas of Lemma \ref{ts5.2.4}
 or the Remark \ref{ts5.2.6}, we obtain:

\sm
\subsubsection{\label{ts5.3.2}Theorem} 
{\it
Let $(N',p)$ and $(N'',q)$ be
two germs of semi--simple  Frobenius 
manifolds with tame 
base--points, Euler fields and flat identities which satisfy
$u^{\prime i}(p) + u^{\prime\prime  j}(q) \neq  
u^{\prime k}(p) + u^{\prime\prime  l}(q)$ 
for $i \neq k$ and $j\neq l$ and 
let the corresponding special initial conditions be given
by $({\Cal V}',\eta')$ and  $({\Cal V}'', \eta'')$, then
the special initial conditions for the Schlesinger
equations corresponding to the tensor product with the flat identity 
and the Euler field of the product chosen as in 
Theorem \ref{ts3.2.1} are given by:
\begin{align}
\eta_{ij} &= \eta'_i \eta''_j\nn\\
v_{ij,kl} &= \delta_{jl} v'_{ik} + \delta_{ik} v''_{jl}.\qed \label{t5.20}
\end{align}
}

\med

\subsubsection{Corollary}
{\it
In the neighborhood of a pair of tame semi--simple base points
the tensor product can be locally given in terms of special initial conditions.
}

\sm

{\bf Proof.}
Since the condition $u^{\prime i}(p) + u^{\prime\prime  j}(q) \neq  
u^{\prime k}(p) + u^{\prime\prime  l}(q)$ 
for $i \neq k$ and $j\neq l$ is an open condition we can always
find a pair of tame semi--simple points satisfying the equation.
Since the tensor product is locally unique up to isomorphism
the Corollary follows.

\med

\subsubsection{Remark}
The virtue of the Corollary above (together with the existence theorem
of the last section and the Theorem on the existence of an Euler field) 
is that it is thus possible to
consider special initial conditions to identify a tensor product
which was originally defined by the tensor product of two germs with
nilpotent base--point.

\sm

In many examples this is exactly the case. For instance in
quantum cohomology as well as in  Saito's unfolding spaces [S, M3]
and the constructions in [D2, DZh]. 

\med

\subsection{Example: 
Special initial conditions for ${\Bbb P}^n \times {\Bbb P}^m$}

Using the Theorem we can calculate the special initial
conditions for ${\Bbb P}^n \times {\Bbb P}^m$ using the results of [MM].
Set $\zeta_n = \exp (\frac{2\pi i}{n+1})$.

\med
\subsubsection{Proposition}
{\it
The point $(x^{00},x^{10},x^{01}, 0, \ldots )$ has canonical 
coordinates $u_{ij}= x^{00} + \zeta_n^i(n+1)e^{\frac{x^{10}}{n+1}}
+\zeta_m^i(m+1)e^{\frac{x^{01}}{m+1}}$

\sm
The special initial conditions at this point corresponding 
to $H_{quant}({\Bbb P}^n \times {\Bbb P}^m$) are given by
\begin{equation}
v_{ij,kl} = -(\frac {\zeta_n^{i-k}}{1-\zeta_n^{i-k}}\delta_{jl}
+\frac {\zeta_m^{j-l}}{1-\zeta_m^{j-l}}\delta_{ik})
\end{equation}
and 
\begin{equation}
\eta_{ij} = \frac{\zeta_n^i\zeta_m^j}{(n+1)(m+1)}
e^{-x^{10}\frac{n}{n+1}-x^{01}\frac{m}{m+1}}
\end{equation}
}

\end {document}